\documentclass[11pt]{article}%
\usepackage{makeidx}
\usepackage{amssymb}
\usepackage[all]{xy}
\usepackage[latin1]{inputenc}
\usepackage{amsfonts}
\usepackage{amsmath}
\usepackage{amssymb}
\usepackage{url}
\usepackage{hyperref}
\usepackage{graphicx}%
\usepackage{mathabx}

\usepackage{amsthm}
\usepackage{url}
\usepackage{hyperref}
\providecommand{\U}[1]{\protect\rule{.1in}{.1in}}

\usepackage{xcolor}

\newtheorem{prop}{Proposition}[section]
\newtheorem{cor}[prop]{Corollary}

\newtheorem{lem}[prop]{Lemma}

\newtheorem{theo}[prop]{Theorem}

\newcommand{\EE}{\mathbb{E}}

\newcommand{\LL}{\mathbb{L}}

\newcommand{\PP}{\mathbb{P}}

\newcommand{\RR}{\mathbb{R}}

\newcommand{\Aa}{ {\cal A }}

\newcommand{\Na}{ {\cal N }}

\newcommand{\Ea}{ {\cal E }}

\newcommand{\Ra}{ {\cal R }}
\newcommand{\Va}{ {\cal V }}
\newcommand{\Ua}{ {\cal U }}
\newcommand{\Fa}{ {\cal F }}

\newcommand{\Qa}{ {\cal Q }}

\newcommand{\Xa}{ {\cal X }}
\newcommand{\Ma}{ {\cal M }}

\newcommand{\Ha}{ {\cal H }}

\newcommand{\Za}{ {\cal Z }}
\newcommand{\Ya}{ {\cal Y }}
\newcommand{\Wa}{ {\cal W }}

\newcommand{\point}{\mbox{\LARGE .}}

\newcommand{\cqfd}{\hfill\blbx \\}
\def\blbx{\hbox{\vrule height 5pt width 5pt depth 0pt}\medskip}

\def \PP{\mathbb{P}}
\def \RR{\mathbb{R}}
\def \EE{\mathbb{E}}

\def \LL{\mathbb{L}}

\def\tr{\mbox{\rm tr}}  

\begin{document}

  \title{On the Stability and the Concentration of Extended Kalman-Bucy filters }
  \author{P. Del Moral, A. Kurtzmann, J. Tugaut}


\maketitle

\begin{abstract}

The exponential stability and the concentration properties of a class of extended Kalman-Bucy filters are analyzed.
New estimation concentration inequalities around partially observed signals are derived in terms of the stability properties of the filters. These non asymptotic exponential inequalities allow to design confidence interval type estimates in terms of the filter forgetting properties with respect to erroneous initial conditions. For  uniformly stable signals,
we also provide explicit non-asymptotic estimates for the exponential forgetting rate of the filters and the associated stochastic Riccati equations w.r.t. Frobenius norms. These non asymptotic exponential concentration and quantitative stability estimates seem to be the first results of this type for this class of 
nonlinear filters.
Our techniques 
combine $\chi$-square concentration inequalities and Laplace estimates with spectral and random matrices  theory, and
 the non asymptotic stability theory of quadratic type stochastic processes.\\



\end{abstract}


\section{Introduction}

The linear-Gaussian stochastic filtering problem has been solved in the beginning of the 1960s by Kalman and Bucy in their seminal articles~\cite{bucy,bucy-2,kalman-60}.  Since this period, Kalman-Bucy filters have became one of the most powerful estimation
algorithm in applied probability, statistical inference, information theory and engineering sciences. The Kalman-Bucy filter is designed to estimate in an optimal way (minimum variance) the internal states of linear-Gaussian time series from a sequence of partial and noisy
measurements. The range of applications goes from tracking, navigation and control to computer vision, econometrics, statistics, finance, and many others. 
For linear-Gaussian filtering problems, the conditional distribution of the internal states of the signal given
the observations up to a give time horizon are Gaussian. The Kalman-Bucy filters and the associated Riccati equation 
coincide with the evolution of the conditional averages and the conditional covariances error matrices of these conditions Gaussian distributions.

Using natural local linearization techniques Kalman-Bucy filters are also currently used to solve
nonlinear and/or non Gaussian signal observation filtering problems. 
The resulting  Extended Kalman-Bucy filter ({\em abbreviated EKF}) often yields powerful and computational efficient 
estimators. Nevertheless it is well known that it fails to be optimal with respect to the minimum variance criteria.
For a more thorough discussion on the origins and the applications of these observer type filtering techniques we refer to  the articles~\cite{luenberger,sonnemann,song} and the
book by D. Simon~\cite{simon}.

 There is a vast literature on the applications and the performance of extended Kalman filter, most on discrete time filtering problems, but very few on the stability properties, none on the exponential concentration properties.
 
In the last two decades, the convergence properties of the EKF have been mainly developed into three different 
but somehow related directions:

 The 
first commonly used approach is to analyze the long-time behaviour of the estimation error between the filter and the partially observed signal. To bypass the fluctuations induced by the signal noise and the observation perturbations,
one natural strategy is to design judicious deterministic observers as the asymptotic
limit of the EKF when the observation and the sensor noise tend to zero. As underlined in~\cite{bonnabel}, in deterministic setting the original covariance matrices of the stochastic signal and the one of the observation perturbations are interpreted as design/tuning type parameters associated with the confidence type matrices of the trusted model and the confidence matrix of associated with the measurements.

For a more detailed discussion on deterministic type observers
as the limit of filters when the sensor and the observation noise tend to zero we refer the reader to the seminal article~\cite{Bensoussan} and the more recent study~\cite{bonnabel}.
Several articles proposed a series of observability and controllability conditions under which the  estimation error of the corresponding discrete time observer converges to zero
~\cite{Bensoussan,boutayed,sonnemann,song}. 
These regularity conditions allow to control the maximal and the minimal eigenvalues of the solution of the Riccati equations (and its inverse). 

One of the drawbacks of this approach is that it gives no precise information on the stochastic EKF but on the limiting noise free-type deterministic observer. On the other hand, up to our knowledge there does not exist any uniform result that allow to quantify the difference between the filter and its asymptotic limit with respect to the time parameter. Another drawback is that
 the initial estimation errors need to be rather small and the signal model close to linear. 

In general practical and stochastic situations, mean square errors do not converge to zero as the time parameter tends to $\infty$. The reasons are two folds: Firstly, the observation noise  of the sensors cannot be totally cancelled. On the other hand the internal signal states are usually only partially observed, and some components may not be fully observable.

A second closely related strategy is to design a Lyapunov function to ensure the stochastic stability of the EKF. Here again these Lyapunov functions are expressed in terms of the inverse of the Riccati equation. These stability properties ensure that the mean square 
estimation error is uniformly bounded w.r.t. the time horizon~\cite{Bensoussan,krener,reif-99,reif-00}.  The regularity conditions are also based on a series of local observability and controllability conditions. As any variance type estimate, these mean square error control are 
somehow difficult to use in practical situations with rather crude confidence interval estimates.

The third and more recent approach is based on the contraction theory developed by W. Lohmiller and J.J.E. Slotine in the seminal articles~\cite{lohmiller-slotine,lohmiller-slotine-2}, and further developed in~\cite{bonnabel}. This approach is also designed to study deterministic type observers.
The idea is  to control the estimation error between a couple of close 
EKF trajectories in a given region w.r.t. the metric induced by the quadratic form associated with the inverse of the solution of the Riccati equation. This approach considers  the partially observed signal as a deterministic system and requires the filter to start in a basin of attraction of the true state.
In summary, these  techniques show that the observer induced by the EKF converges locally exponentially to the state of the signal
when the quadratic form induced by the inverse of the Riccati equation is sufficiently regular and under appropriate observability and controllability conditions.

The objective of this article is to complement these three approaches with a  
 novel stochastic analysis
based on exponential concentration inequalities and uniform $\chi$-square type estimates 
for stochastic quadratic type processes. 

Our regularity conditions are somehow stronger than the ones discussed in the above referenced articles
but they do not rely on suitable local initial conditions nearby the true signal state. Last but not least our methodology applies to stochastic filtering problems, not to deterministic type observers. 

In our framework the signal process is a uniformly and 
 exponentially stable Langevin type diffusion, and the sensor function is the identity matrix up to a change of basis. 
 
 In this apparently simple nonlinear filtering problem the quantitative analysis of the EKF exponential stability is based on sophisticated probabilistic tools.
The complexity of these stochastic processes can be measured by the fact that 
the  EKF is a nonlinear diffusion process equipped with a diffusion correlation matrix satisfying a coupled nonlinear and 
stochastic Riccati equation.

 This study has been motivated by one of our recent research project on the refined convergence analysis Ensemble type
Kalman-Bucy filters. To derive some useful uniform convergence results with respect to the time horizon we have shown in~\cite{dt-1-2016} that the signal process needs to be uniformly stable and fully observed by some noisy sensor. These rather strong conditions cannot be relaxed even for linear Gaussian filtering models. We plan to extend these results for nonlinear filtering models based on the non asymptotic estimates presented in this article.

In this context we present new exponential concentration inequalities
 to quantify the stochastic stability of the EKF. 
 They allow to derive confidence intervals for 
 the deviations of the stochastic flow of the EKF around the internal states of the partially observed signal. These  estimates also show that
the fluctuations induced by any erroneous initial condition tend to zero as the time horizon tends to $+\infty$.

Our second objective is to develop a non asymptotic quantitative analysis of the stability properties of the EFK.
In contrast to the linear-Gaussian case discussed in~\cite{dt-1-2016}, the Riccati equation associated with the EFK depends on the states of the filter. The resulting system is a nonlinear stochastic process evolving in multidimensional
inner product spaces. To analyze these complex models we develop a stability theory of quadratic type stochastic processes. Our main contribution is a non asymptotic $\LL_p$-exponential stability theorem. This theorem shows that the $\LL_p$-distance between
two solutions of the EKF and the stochastic Riccati equation with possibly different initial conditions converge to zero as the time horizon tends to $+\infty$.  We also provide a non asymptotic estimate of the exponential decay rate.

The rest of the article is organized as follows: 

In the next two sections,
Section~\ref{sec-description} and Section~\ref{statement-section}, we present the nonlinear filtering models 
discussed in the article and we state the main results developed in this work. Section~\ref{stab-chi-2} is concerned
with the stability properties of quadratic type processes. This section presents the main technical results used in the further development of the article. Most of the technical proofs are provided in the appendix. Section~\ref{stab-prop-section}
is dedicated to the stochastic stability properties of the signal and the EKF. The end of the article is mainly concerned with the proofs of the two main theorems presented in Section~\ref{statement-section}.

\subsection{Description of the models}\label{sec-description}

This section presents the nonlinear filtering models in this article. We also discuss and illustrate our regularity conditions
with several classes of Langevin type signal processes partially observed by  noisy sensors.

Consider a time homogeneous nonlinear filtering problem of the following form
\begin{equation}\label{filtering-model-ref}
\left\{
\begin{array}{rcl}
dX_t&=&A(X_t)~dt~+~R^{1/2}_{1}~dW_t\\
dY_t&=&BX_t~dt~+~R_2^{1/2}~dV_t
\end{array}
\right.\quad\mbox{\rm and we set $\Fa_t=\sigma\left(Y_s,~s\leq t\right)$. }
\end{equation}
In the above display,  $(W_t,V_t)$ is an $(r_1+r_2)$-dimensional Brownian motion, $X_0$ is a $r_1$-valued
Gaussian random vector with mean and covariance matrix $(\EE(X_0),P_0)$ (independent of $(W_t,V_t)$), the symmetric matrices $R^{1/2}_{1}$ and $R^{1/2}_{2}$ are invertible, $B$ is an  $(r_2\times r_1)$-matrix, and $Y_0=0$. The 
drift of the signal 
 is a differentiable vector valued function $
A~:~x\in \RR^{r_1}\mapsto A(x)\in \RR^{r_1}
$ with a Jacobian denoted by $\partial A~:~x\in \RR^{r_1}\mapsto A(x)\in \RR^{(r_1\times r_1)}$. 

The Extended Kalman-Bucy filter associated with the filtering problem (\ref{filtering-model-ref}) is defined by the evolution equations
\begin{eqnarray}
d\widehat{X}_t&=&A(\widehat{X}_t)~dt+P_tB^{\prime}R^{-1}_{2}~\left[dY_t-B\widehat{X}_tdt\right]\quad\mbox{\rm with}\quad 
\widehat{X}_0=\EE(X_0)\nonumber
\\
\partial_tP_t&=&\partial A(\widehat{X}_t)P_t+P_t\partial A(\widehat{X}_t)^{\prime}+R_1-P_t~SP_t\quad\mbox{\rm with}\quad S:=B^{\prime}R^{-1}_2B
\label{EKF}
\end{eqnarray}
where $B^{\prime}$ stands for the transpose of the matrix $B$.
For nonlinear signal processes the random matrices $P_t$ cannot be interpreted as the error covariance matrices.  
Nevertheless, rewriting the EKF in terms of the signal process we have
\begin{eqnarray*}
d(X_t-\widehat{X}_t)&=&[(A(X_t)-A(\widehat{X}_t))-P_tS(X_t-\widehat{X}_t)]~dt+R^{-1/2}_{1}dW_t
-P_tB^{\prime}R^{-1/2}_{2}dV_t
\end{eqnarray*}
Replacing $(A(X_t)-A(\widehat{X}_t))$  by the first order approximation
$\partial A(\widehat{X}_t)(X_t-\widehat{X}_t)$ we define a process
$$
d\widetilde{X}_t:=[\partial A(\widehat{X}_t)-P_tS]\widetilde{X}_t~dt+R^{-1/2}_{1}dW_t
-P_tB^{\prime}R^{-1/2}_{2}dV_t.
$$
It is a simple exercise to check that the solution of the Riccati equation (\ref{EKF}) coincides with the $\Fa_t$-conditional
covariance matrices of $\widetilde{X}_t$; that is, for any $t\geq 0$ we have 
$
P_t=\EE\left(\widetilde{X}_t\widetilde{X}_t^{\prime}~|~\Fa_t\right)
$.

\subsubsection{Langevin-type signal processes}
 In the further development of the article we assume that the  Jacobian matrix of $A$ satisfies the following regularity conditions:
\begin{equation}\label{regularity-condition-varpiAB}
\left\{
 \begin{array}{rcl}
 -\lambda_{\partial A}&:=&\sup_{x\in \RR^{r_1}} \rho(\partial A(x)+\partial A(x)^{\prime})<0\\
 \\
\Vert \partial A(x)-\partial A(y)\Vert&\leq &\kappa_{\partial A}~\Vert x-y\Vert
\quad\mbox{\rm
 for some $\kappa_{\partial A}<\infty$.}
 \end{array}\right.
\end{equation}
where $\rho(P):=\lambda_{\tiny max}(P)$ stands for the maximal eigenvalue of a symmetric matrix $P$. In the above display $\Vert \partial A(x)-\partial A(y)\Vert$ stands for the $\LL_2$-norm of the matrix operator $(\partial A(x)-\partial A(y))$, and $\Vert x-y\Vert$ the Euclidean distance between $x$ and $y$.
A Taylor first order expansion shows that 
\begin{equation}\label{regularity-conditon}
(\ref{regularity-condition-varpiAB})\Longrightarrow \langle x-y,A(x)-A(y)\rangle\leq -\lambda_A~\Vert x-y\Vert^2\quad\mbox{\rm with $\lambda_A\geq \lambda_{\partial A}/2> 0$. }
\end{equation} 

The above rather strong conditions ensure the contraction needed to ensure the stability of the EFK. They are also used to derive uniform estimates 
w.r.t. the time horizon for Ensemble Kalman-Bucy particle filters~\cite{dtk-2016}. For linear systems $A(x)=Ax$, associated with some matrix $A$, the parameters $\lambda_{A}=\lambda_{\partial A}/2$ coincide with the logarithmic norm of $A$.

The prototype of signals satisfying these conditions are multidimensional diffusions
with  drift functions $(A,\partial A)=(-\partial \Va,-\partial^2 \Va)$ associated with a gradient Lipschitz  strongly convex confining potential $ \Va~:~x\in \RR^{r_1}\mapsto \Va(x)\in [0,\infty[$. The logarithmic norm condition (\ref{regularity-condition-varpiAB}) is met as soon as $\partial^2 \Va\geq v~Id$ with $v=2\vert\lambda_{\partial A}\vert$.   Equivalently the smallest eigenvalue $\lambda_{\tiny min}(\partial^2 \Va(x))$ of the Hessian  is uniformly lower bounded by $v$.
In this case  (\ref{regularity-condition-varpiAB}) is met with $\lambda_{\partial A}=v/2$.  These conditions are fairly standard in the stability theory of nonlinear diffusions, we refer the reader to the review article~\cite{villani}, and the references therein. Choosing 
 $R_1=\sigma_1^2~Id$ and $A=-\beta\partial  \Va$, for some $\beta,\sigma_1\geq 0$ the signal process $X_t$ resumes to a multidimensional
 Langevin-diffusion
\begin{equation}\label{langevin-example}
dX_t=-\beta~ \partial  \Va(X_t)~dt+\sigma_1~dW_t.
\end{equation}
This process is reversible w.r.t. the invariant distribution $\mu_\beta$, where $\mu_\beta$ is the  probability distribution
on $\RR^{r_1}$ given by
$$
\mu_{\beta}(dx)=\frac{1}{\Za_{\beta}}~\exp{\left(-\frac{2\beta}{\sigma_1^2} \Va(x)\right)}~dx\quad \mbox{\rm with}\quad \Za_{\beta}=\int~\exp{\left(-\frac{2\beta}{\sigma_1^2} \Va(x)\right)}~dx\in ]0,\infty[.
$$
In the above display $dx$ stands for the Lebesgue measure on $\RR^{r_1}$. 
The Lipschitz-continuity condition of the Hessian $\partial^2 \Va$ introduced in (\ref{regularity-condition-varpiAB}) ensures the continuity of the stochastic Riccati equation (\ref{EKF}) w.r.t. the fluctuations around the random states $\widehat{X}_t$. We illustrate this condition with a
nonlinear example given by the function $$
 \Va(x)=\frac{1}{2}~\langle \Qa_1 x, x\rangle+ \langle q,x\rangle+\frac{1}{3}~\langle \Qa_2 x,x\rangle^{3/2}
$$ with some symmetric positive definite matrices $(\Qa_1,\Qa_2)$ and some given
vector $q\in\RR^{r_1}$. In this case we have
\begin{eqnarray*}
\partial  \Va(x)&=&q+\Qa_1x+\langle \Qa_2 x,x\rangle^{1/2}~\Qa_2 x,\\
\partial^2  \Va(x)&=&\Qa_1+\langle \Qa_2 x,x\rangle^{1/2}~\Qa_2+\langle \Qa_2 x,x\rangle^{-1/2}~\Qa_2 xx^{\prime}\Qa_2\geq \Qa_1.
\end{eqnarray*}
In this situation we have
\begin{equation}\label{example-non-linear}
\Vert \partial^2  \Va(x)-\partial^2  \Va(y)\Vert\leq 2~\Vert \Qa_2\Vert^{3/2}~\Vert y-x\Vert.
\end{equation}
This shows that conditions (\ref{regularity-condition-varpiAB}) are met with the parameters 
$$(\lambda_{\partial A},\kappa_{\partial A})=\beta~\left(2^{-1}\lambda_{\tiny min}(\Qa_1),~2\lambda_{\tiny max}^{3/2}(\Qa_2)\right).$$ A proof of (\ref{example-non-linear}) is provided in the appendix on page~\pageref{proof-formula-Vnl}.
More generally these regularity conditions also hold if we replace in (\ref{langevin-example}) the parameter $\sigma_1$ by 
any choice of covariance matrice $R_1$. Also observe that the Langevin diffusion associated with the null form $\Qa=0$ coincides with the conventional linear-Gaussian filtering problem discussed in~\cite{dt-1-2016}. Stochastic gradient-flow diffusions of the form
(\ref{langevin-example}) arise in a variety of application domains. In mathematical finance and mean field game theory~\cite{carmona-fouque,fouque-sun}, these Langevin models describe the interacting-collective behaviour of $r_1$-individuals. For instance in the Langevin model 
discussed in~\cite{fouque-sun} the state variables $X_t=\left(X^i_t\right)_{1\leq i\leq r_1}$ represent the log-monetary reserves of $r_1$ banks
lending and borrowing to each other.  The quadratic potential function is given by
$$
\langle \Qa_1 x, x\rangle=\sum_{1\leq i\leq r_1}\left(x_i-\frac{1}{r_1}\sum_{1\leq j\leq r_1}~x_j\right)^2\Rightarrow~\Qa_1
\succ \left(1-\frac{1}{r_1}\right)~I_{r_1}.
$$
In this context, the parameter $\beta$ represents the mean-reversion rate between banks. More general interacting potential functions can be considered.  Mean field type diffusion processes are also used to design low-representation of fluid flow velocity fields. These 
vortex-type particle filtering problems are developed in some details in the pionnering articles by E. M\'emin and his co-authors~\cite{cuzol1,cuzol2,cuzol3,papadakis}. These probabilistic interpretations of the 2d-incompressible Navier-Stokes equation 
represent the vorticity map as a mixture of basis functions centered around
each vortex.  

In this connexion, we mention that our approach also applies to interacting diffusion  gradient flows described 
by a potential function of the form
$$
\Va(x)=\sum_{1\leq i\leq r_1}~\Ua_1(x_i)+\sum_{1\leq i\not=j\leq r_1}~\Ua_2(x_i,x_j)
$$
for some gradient Lipschitz  strongly convex confining potential $ \Ua_i:\RR^i\mapsto [0,\infty[$, $i=1,2$. In this situation, we have
\begin{equation}\label{U1U2U}
\partial^2 \Ua_1\geq u_1\quad\mbox{\rm and}
\quad \partial^2\Ua_2\geq u_2~I_2\Longrightarrow~\partial^2\Va\succ v~I_{r_{1}}\quad \mbox{\rm with}\quad v:=(u_1+(r_1-1)u_2)>0.
\end{equation}
We further assume that
\begin{eqnarray*}
\vert \partial^2\Ua_1(x_1)-\partial^2\Ua_1(y_1)\vert&\leq& \kappa_{\partial^2\Ua_1}~\vert x_1-x_2\vert,\\
\Vert \partial^2\Ua_2(x_1,x_2)-\partial^2\Ua_2(y_1,y_2)\Vert &\leq &\kappa_{\partial^2\Ua_2}~\Vert (x_1,x_2)-(y_1,y_2)\Vert .
\end{eqnarray*}
In this case,  we have
\begin{equation}\label{2nd-formula-Lipschitz}
\Vert \partial^2\Va(x)- \partial^2\Va(y)\Vert\leq \kappa_{\partial^2\Va}~\Vert x-y\Vert\quad\mbox{\rm with}\quad
 \kappa_{\partial^2\Va}:=\kappa_{\partial^2\Ua_1}+\kappa_{\partial^2\Ua_2}~(r_1-1)~ \sqrt{2(r_1-1)}.
\end{equation}
This shows that conditions (\ref{regularity-condition-varpiAB}) are met with $$(\lambda_{\partial A},\kappa_{\partial A})=\beta~\left(2^{-1}(u_1+(r_1-1)u_2),~\kappa_{\partial^2\Ua_1}+\kappa_{\partial^2\Ua_2}~(r_1-1)~ \sqrt{2(r_1-1)}\right).$$ 
The detailed proofs of (\ref{U1U2U})-(\ref{2nd-formula-Lipschitz}) are provided in the appendix on page~\pageref{proof-formula-Vnl}.

\subsubsection{Observability conditions}
When the observation variables are the same as the ones of the signal; the signal observation has the same dimension as the signal and resumes to
 some equation of the form
\begin{equation}\label{fully-observed-sensor}
dY_t=b~ X_t~dt+\sigma_2~dV_t
\end{equation}
for some parameters $b\in\RR$ and $\sigma_2\geq 0$.  These sensors are used in data grid-type assimilation
problems when measurements can be evaluated at each cell. These fully observed models are discussed  in~\cite[Section 4]{harlim-hunt} in the context of the Lorentz-96 filtering problems. These observation processes are also used in the article~\cite{berry-harlim} for application to nonlinear and multi-scale filtering problem. In this context, the observed variables represent the slow components of the signal.

 For partially observed signals we cannot expect any stability properties of the EKF and the EnKF without introducing some structural conditions of observability and controllability on the signal-observation equation (\ref{filtering-model-ref}). Observe that the EKF  equation (\ref{EKF}) implies that
 $$
d(\widehat{X}_t-X_t)=\left[(A(\widehat{X}_t)-A(X_t))-P_tS(\widehat{X}_t-X_t)\right]~dt+P_{t}~C^{\prime}R^{-1/2}_{2}~dV_t+R^{1/2}_{1}~dW_t.
 $$
 This equation shows that the stability properties of this process depends on the nature of the real eigenvalues of the symmetric 
 matrices $(A(x)-PS)_{\tiny sym}$, with $x\in\RR^{r_1}$. In contrast with the conventional Kalman-Bucy filter the Riccati equation
 (\ref{EKF}) is a stochastic equation. As a result, the stability property of the EKF is not
induced by some kind of observability condition that ensures the existence of a steady state deterministic covariance matrix. 

The random fluctuations of the matrices $\partial A(\widehat{X}_t)$ entering in the Riccati equation  (\ref{EKF}) may corrupt the stability in the EKF, even if the linearized filtering problem around some chosen state is observable and controllable. For a more thorough discussion on the stability properties of Kalman-Bucy filters and Riccati equations for linear Gaussian filtering problems we refer the reader to~\cite{abou,anderson,freiling,sontag,wonham,zelikin}.

The stability analysis of diffusion processes is always much more documented than the ones on their possible divergence.
For instance, in contrast with conventional Kalman-Bucy filters, the stability properties of the EKF are not
induced by some kind of observability or controllability condition. The only known results in this direction is the recent pioneering work by X. T. Tong, A. J. Majda and D. Kelly~\cite{tong}  in the context of discrete generation
Ensemble Kalman filters. One of the main assumptions of the article is that the sensor-matrix 
has full rank. The authors also provide a concrete numerical example of filtering problem with sparse observations for which the EnKF experiences a catastrophic divergence.

In a recent article~\cite{dt-1-2016}, in the context of linear drift functions we also show that the uniform propagation of chaos properties of EnKF require 
strong signal stability properties and the same type of observability conditions. 

They are some strong similarities between the EKF and
the EnKF:

 The first one comes from the fact that the predictable part of the EKF is stochastic and nonlinear. The predictable part of the 
EnKF also depend on stochastic covariance matrices. These interaction functions are clearly nonlinear in the internal states of the particle system.

The second one comes from the fact that the Riccati equation associated with the EKF is stochastic. The stochastic perturbation theorem~obtained in~\cite[Theorem 3.1]{dt-1-2016} also shows that the sample covariance matrices satisfy a stochastic diffusion type Riccati equation. 

Without any strong observability conditions, these stochastic nonlinearities may corrupt severally the stability of the EKF.
In the further development of the article we shall assume that the sensor function has the same form as the one discussed in~\cite{dt-1-2016}. 
More precisely, we assume that the following observability condition is satisfied:
\begin{equation}\label{new-condition}
\hskip-5cm(\mbox{\rm S})\hskip3cm \qquad S=\rho(S)~Id\quad\mbox{\rm for some}\quad \rho(S)>0.
\end{equation}

The fully observed model discussed in (\ref{fully-observed-sensor}) clearly satisfies
condition (\ref{new-condition})  with the parameter $\rho(S)=(b/\sigma_2)^2$.  
As mentioned above, in the context of linear-Gaussian filtering problems
this condition is also essential to ensure the uniform convergence of Ensemble Kalman-Bucy filter w.r.t. the 
time parameter.
Section 4 in the article~\cite{dt-1-2016} provides a detailed discussion on spectral estimates and 
semigroup contraction inequalities based on this condition. A
geometric description of global divergence regions in the set of positive covariances matrices is also provided in the context of $2$-dimensional
partially observed filtering problems.

Last but not least, we mention that (\ref{new-condition}) is satisfied when the filtering problem is similar to the ones discussed above; that is, up to a change of basis functions.
For instance (\ref{new-condition}) is met with $S=I_{r_1}$ for sensors with orthonormal  matrices $BR^{-1/2}_2$.  Under this condition, up to a change of observation basis, the observation process reduces to 
$$
\overline{Y}_t:=B^{\prime}R^{-1}_2Y_t\Rightarrow  d\overline{Y}_t=\rho(S)~X_t~dt+~d\overline{V}_t
$$
with an $r_1$-dimensional Wiener process
$
\overline{V}_t:=B^{\prime}R^{-1/2}_2~V_t
$. Inversely, any filtering problem (\ref{filtering-model-ref}) with $r_1=r_2$ and s.t.  $(R_2^{-1/2}B)$ is invertible can be turned into that form.
To check this claim we observe that
$$
\Ya_t:=R^{-1/2}_2Y_t\quad\mbox{\rm and}\quad \Xa_t:=R^{-1/2}_2BX_t
\Longrightarrow \left\{\begin{array}{rcl}
d\Xa_t&=&\Aa(\Xa_t)~dt+\Ra_1^{1/2}~dW_t\\
d\Ya_t&=&\Xa_t~dt+dV_t
\end{array}\right.
$$
with the drift function
$$
\Aa:=(R^{-1/2}_2B)\circ A\circ (R^{-1/2}_2B)^{-1}\quad \mbox{\rm and the matrix}\quad
\Ra_1:=R^{-1/2}_2 BR_1B^{\prime}R^{-1/2}_2.
$$
In this situation the filtering model $(\Xa_t,\Ya_t)$ satisfies (\ref{new-condition}) with $\rho(S)=1$. In addition, we have
$$
A=(R^{-1/2}_2B)^{-1}\circ \partial U \circ (R^{-1/2}_2B) \Rightarrow
(\Aa,\partial \Aa)=(\partial U,\partial ^2U).
$$

\subsection{Statement of the main results}\label{statement-section}
We let $\phi_t(x):=X_t$ and $\varphi_t(x):=x_t$ be the stochastic and the deterministic 
flows of the stochastic and the deterministic systems
$$
\left\{
\begin{array}{rcl}
dX_t&=&A(X_t)~dt+R^{1/2}_1dW_t\\
\partial_tx_t&=&A(x_t)\hskip3cm
\mbox{\rm
starting at $x_0=\varphi_0(x)=x=X_0=\phi_0(x)$.}
\end{array}\right.
$$
We also let $\overline{\Phi}_t :=(\Phi_t,\Psi_t)$ be the stochastic flow associated with the EKF and the Riccati
stochastic differential equations; that is
$$
\overline{\Phi}_t(\widehat{X}_0,P_0)=\left(\Phi_t(\widehat{X}_0,P_0),\Psi_t(\widehat{X}_0,P_0)\right):=\left(\widehat{X}_t,P_t\right).
$$

Given  $(r_1\times r_2)$ matrices $P,Q$ we define the Frobenius inner product
$$
\langle P,Q\rangle=\tr(P^{\prime}Q)\quad\mbox{\rm and the associated norm}\quad
\Vert P\Vert_F^2=\tr(P^{\prime}P)
$$
where $\tr(C)$ stands for the trace of the matrix $C$. 
We also equip the product space $\RR^{r_1}\times \RR^{r_1\times r_1}$ with the inner product 
$$
\langle (x_1,P_1),(x_2,P_2)\rangle:=\langle x_1,x_2\rangle+\langle P_1,P_2\rangle\quad
\mbox{\rm and the norm}\quad \Vert (x,P)\Vert^2:=\langle (x,P),(x,P)\rangle .
$$

We recall the $\chi$-square Laplace estimate
\begin{equation}\label{chi-2}
\EE\left(\exp{\left[ {\Vert X_0-\widehat{X}_0\Vert^{2}}/{\chi(P_0)}\right]}\right)\leq e\qquad\mbox{\rm with}\quad
\chi(P_0):=4r_1\rho(P_0).
\end{equation}
The proof of (\ref{chi-2}) and more refined estimates are housed in the appendix. We have the rather crude almost sure estimate 
\begin{eqnarray*}
\partial_t\tr\left(P_t\right)&=&\tr\left((\partial A(\widehat{X}_t)+\partial A(\widehat{X}_t)^{\prime})P_t\right)+\tr(R_1)-\tr(SP_t^2)\\
&\leq & -\lambda_{\partial A}~\tr\left(P_t\right)+\tr(R_1).
\end{eqnarray*}
This readily yields
the upper bound
\begin{equation}\label{riccati-bound}
\tr(P_t)\leq \tau_t(P):=e^{-\lambda_{\partial A}t}~\tr\left(P_0\right)+{\tr(R_1)}/{\lambda_{\partial A}}\quad\Rightarrow\quad\sup_{t\geq 0}\tr(P_t)\leq \tr\left(P_0\right)+{\tr(R_1)}/{\lambda_{\partial A}}.
\end{equation}

Most of the analysis developed in the article relies on the following quantities:
\begin{equation}\label{def-sigma-pi}
\sigma^2_{\partial A}:=
1+2~\pi_{\partial A}\quad\mbox{\rm with}\quad \pi_{\partial A}(t):= \tau_t^2(P)~\rho(S)~\tr(R_1)^{-1}\longrightarrow_{t\rightarrow\infty} \pi_{\partial A}:=\frac{\rho(S)}{\lambda_{\partial A}}~
\frac{\tr(R_1)}{\lambda_{\partial A}}.
\end{equation}

%
Our first main result concerns the  stochastic stability of the EKF and is described in terms of the function
$$
\delta\in [0,\infty[\mapsto \varpi(\delta):= \frac{e^2}{\sqrt{2}}~\left[\frac{1}{2}+\left(\delta+\sqrt{\delta}\right)\right].
$$

More precisely we have the following
exponential concentration theorem.

\begin{theo}
For any initial states 
$(x,\widehat{x},p)\in\RR^{r_1+r_1+(r_1\times r_1)}$ and any time horizon $t\in [0,\infty[$, and
any $\delta\geq 0$
 the probabilities  of the following
events are greater than $1-e^{-\delta}$:
\begin{eqnarray}
\Vert \phi_t(x)-\varphi_t(x)\Vert^2&\leq& ~\varpi(\delta)~\frac{\tr(R_1)}{\lambda_A},\label{concentration-stab-signal}
\\
&&\nonumber\\
\Vert \phi_t(x)-\Phi_t(\widehat{x},p)\Vert^2
\displaystyle&\leq& 4~\varpi(\delta)~\frac{\tr(R_1)}{\lambda_A}~\sigma^2_{\partial A}
\label{concentration-stab}\\
&&+2
e^{-\lambda_{\partial A}t}~\Vert x-\widehat{x}\Vert^2+8\varpi(\delta)~
\frac{\vert e^{-\lambda_At}- e^{-\lambda_{\partial A}t}\vert}{\vert \lambda_A-\lambda_{\partial A}\vert}~
\rho(S)~\tr(p)^2 .
\nonumber
\end{eqnarray}
\end{theo}

The proofs of the concentration inequalities (\ref{concentration-stab-signal}) and (\ref{concentration-stab}) are 
provided respectively in Section~\ref{signal-sec} and Section~\ref{lp-section}. See also 
Theorem~\ref{prop-laplace-signal-dyns} and Theorem~\ref{theo-XwX}
for related Laplace $\chi$-square estimates of time average distances.

The role of each quantity in (\ref{concentration-stab-signal}) and (\ref{concentration-stab}) is clear. The size of the ``confidence events" are proportional to the signal or the observation perturbations, and inversely proportional to the stability rate of the systems. More interestingly, formula (\ref{concentration-stab}) shows that the impact of the initial conditions is exponentially small when the time horizon increases.

Our next objective is to better understand the stability properties of the EKF and the corresponding stochastic Riccati equation.
 To this end, it is convenient to strengthen our regularity conditions. We further assume that
 \begin{equation}\label{reg-hyp}
 \lambda_{\partial A}>\sqrt{2\kappa_{\partial A}\tr(R_1)}\vee (4\rho(S))
 \end{equation}
 and  for some  $\alpha>1$
\begin{equation}\label{regularity-conditon-HAS-0}
4e\alpha~\sqrt{\frac{\rho(S)}{\lambda_{\partial A}}}~\frac{\tr(R_1)}{\lambda_{A}}~\left[1+2~\frac{\tr(R_1)}{\lambda_{\partial A}}
\frac{\rho(S)}{\lambda_{\partial A}}\right]<1.
\end{equation}

In contrast with the linear-Gaussian case, the Riccati equation (\ref{EKF}) depends on the internal states of the EKF.
As a result its stability properties are characterized by a stochastic Lyapunov exponent that depends on the 
random trajectories of the filter as well as on the signal-observation processes.
Condition (\ref{regularity-conditon-HAS-0}) is a technical condition that allows to control uniformly the fluctuations of these stochastic
exponents with respect to the time horizon.
By (\ref{reg-hyp}) this condition is met as soon as
$$
\alpha e~{\tr(R_1)}~\left[1+\frac{1}{8}~\frac{\tr(R_1)}{\rho(S)}\right]<\lambda_{A}/2.
$$

Loosely speaking, when the signal is not sufficiently stable the erroneous initial conditions of EKF may be too sensitive to small perturbations of the sensor. When the exponential decay to equilibrium of the signal is stronger than these spectral instabilities 
the EKF and the corresponding stochastic Riccati equations are stable and forgets any erroneous initial conditions.

We set
$$
\Delta_t:=\Vert\overline{\Phi}_t\left(\widehat{X}_0,P_0\right)-
\overline{\Phi}_t\left(\widecheck{X}_0,\widecheck{P}_0\right)\Vert^2\quad
$$
and
\begin{eqnarray*}
 \Lambda/\lambda_{\partial A}&:=&
1-2~\frac{\kappa_{\partial A}}{\lambda_{\partial A}}~\frac{\tr(R_1)}{\lambda_{\partial A}}-~
\sqrt{\frac{\rho(S)}{\lambda_{\partial A}}}\left[1-\frac{3}{4}
\sqrt{\frac{\rho(S)}{\lambda_{\partial A}}}\right]\\
&\geq &\frac{1}{2}-2~\frac{\kappa_{\partial A}}{\lambda_{\partial A}}~\frac{\tr(R_1)}{\lambda_{\partial A}}>0.
\end{eqnarray*}

We are now in position to state our second main result.

\begin{theo}\label{stab-theo}
When $\lambda_{\partial A}>0$ we have the uniform estimates
$$
\forall n\geq 1\qquad \sup_{t\geq 0}\EE(\Delta_t^n)<\infty.
$$
Assume conditions (\ref{reg-hyp}) and  (\ref{regularity-conditon-HAS-0}) are satisfied for  some $\alpha> 1$. In this situation,
for any $\epsilon\in ]0,1]$  there exists some time horizon $s$ such that
for any $t\geq s$ we have the almost sure contraction estimate 
$$
 \delta:=\frac{1}{2}
\sqrt{\frac{\lambda_{\partial A}}{\rho(S)}}~(~> 1)~\Longrightarrow~
\EE\left(\Delta_t^{\delta/2}~|~\Fa_s\right)^{2/\delta}
\displaystyle\leq  \Za_{s}~\displaystyle \exp{\left[-\left(1-\epsilon\right)\Lambda(t-s)\right]}~\Delta_s
$$
for some random process $\Za_{t}$ s.t.
$
\sup_{t\geq 0}{\EE\left( 
\Za^{\alpha\delta}_{t}
\right)}<\infty
$.

 \end{theo}

Theorem~\ref{stab-theo} readily implies 
  the stability  $\delta$-moment Lyapunov exponent estimates
$$
\liminf_{t\rightarrow\infty}-\frac{1}{\delta t}\log{\EE(\Delta_t^{\delta})}\geq \Lambda .
$$
In addition we have the non asymptotic estimates
$$
\EE\left(\Delta_t^{\delta/2}\right)^{2/\delta}\leq~\nu(\alpha)~\displaystyle \exp{\left\{-\left(1-\epsilon\right)\Lambda~(t-s)\right\}}~\EE\left[
\Delta_{s}^{
\delta\alpha/(2\alpha-1)}\right]^{(2\alpha-1)/(\delta\alpha)}
$$
with
$$
\nu(\alpha):=\sup_{t\geq 0}{\EE\left( 
\Za_t^{\alpha\delta}
\right)^{{1}/{(\delta\alpha)}}}<\infty .
$$

We end this section with some comments on our regularity conditions.
Notice that $\Lambda$ does not depend on the parameter $\delta$ nor on $\rho(S)$.
As mentioned above, we believe that these technical conditions can somehow be relaxed. These conditions are stronger than the ones discussed in~\cite{dt-1-2016} for linear-Gaussian models. 
In contrast with the linear case, the Riccati equation in nonlinear settings is a stochastic process in matrix spaces. For this class of models, these technical conditions
are used to control the fluctuations of the stochastic Riccati equation entering into the EKF.

\section{Stability properties of quadratic type processes}\label{stab-chi-2}

Let $(\Ua_t,\Va_t,\Wa_t,\Ya_t)$ be some non-negative processes  defined on 
 some probability space  $(\Omega,\Fa,\PP)$  equipped with a filtration $\Fa=(\Fa_t)_{t\geq 0}$
of $\sigma$-fields. Also let $(\Za_t,\Za^+_t)$  be some processes and $\Ma_t$  be some continuous  $\Fa_t$-martingale. We use the notation 
\begin{equation}\label{stoch-ineq}
d\Ya_t\leq \Za_t^+~dt+d\Ma_t\Longleftrightarrow \left(d\Ya_t=\Za_t~dt+d\Ma_t\quad\mbox{\rm with}\quad \Za_t\leq \Za^+_t\right)
\end{equation}
Let us mention some useful properties of the above stochastic inequalities.

Let $(\overline{\Ya}_t,\overline{\Za}_t^+,\overline{\Za}_t,\overline{\Ma}_t)$ be another collection of processes satisfying
the above inequalities. In this case it is readily checked that
$$
d(\Ya_t+\overline{\Ya}_t)\leq (\Za_t^++\overline{\Za}_t^+)~dt+d(\Ma_t+\overline{\Ma}_t)
$$
and
$$
d(\Ya_t\overline{\Ya}_t)\leq \left[\overline{\Za}_t^+\Ya_t+ \Za_t^+\overline{\Ya}_t+\partial_t\langle \Ma,\overline{\Ma}\rangle_t\right]~dt+\Ya_t~d\overline{\Ma}_t+\overline{\Ya}_t~d\Ma_t .
$$

Let $(\Ha,\langle\point,\point\rangle)$ be some inner product space, and let $\Aa_t~:~x\in\Ha\mapsto \Aa_t(x)\in \Ha$ be a linear operator-valued stochastic process 
with finite logarithmic norm $\rho(\Aa_t)<\infty$.
Consider  an $\Ha$-valued stochastic process $\Xa_t$ such that
\begin{equation}\label{H-general}
d\Vert\Xa_t\Vert^2\leq \left[\langle\Xa_t,\Aa_t\Xa_t\rangle+\Ua_t\right]~dt+d\Ma_t
\end{equation}
for some continuous $\Fa_t$-martingale $\Ma_t$ with angle bracket satisfying the following property
$$
\partial_t \langle\Ma\rangle_t\leq \Va_t~\Vert\Xa_t\Vert^2+\Wa_t~\Vert\Xa_t\Vert^4 .
$$
This section is concerned with the long-time quantitative behaviour of the above quadratic type processes. The main difficulty here comes from the fact that $\Aa_t$ is a stochastic flow of operators. As a result we cannot apply conventional Lyapunov techniques based on Dynkin's formula, supermartingale theory and/or more conventional Gronwall type estimates.

Next theorem provides a way to estimate these processes in terms of geometric type  processes and exponential martingales.

\begin{theo}\label{gronwall-2dim}
When $\Ua_t=0=\Va_t$ we have the almost sure estimate
$$
\Vert \Xa_t\Vert^2\leq \Vert \Xa_0\Vert^2~\exp{\left(\int_0^t\rho(\Aa_s)~ds\right)}~\exp{\left(\int_0^t\sqrt{\Wa_s}~d\Na_s-\frac{1}{2}~\int_0^t\Wa_s~ds\right)}
$$
with  a martingale $\Na_t$ s.t. $\partial_t\langle  \Na\rangle_t\leq 1$.
More generally, for any $n\geq 1$ we have
\begin{equation}\label{gronwall-n+}
\begin{array}{l}
\EE\left(\Vert\Xa_t\Vert^{n}~|~\Fa_0\right)^{2/n}
\displaystyle\leq \EE\left[\exp{\left(n~\int_0^t\left\{\rho(\Aa_s)~+\frac{(n-1)}{2}~ \Wa_s\right\}~ds~|~\Fa_0\right)}\right]^{1/n}\\
\\
\hskip3cm\times\displaystyle\left\{
\Vert\Xa_0\Vert^{2}+\int_0^t~\left(
 \EE\left[\overline{\Ua}_s^n~|~\Fa_0\right]^{1/n}+\frac{(n-1)}{2}~\EE\left[\overline{\Va}_s^n~|~\Fa_0\right]^{1/n}\right)~ds
\right\}
\end{array}\end{equation}
with the rescaled processes 
$$
\overline{\Ua}_t/\Ua_t:=\exp{\left(-\int_0^t\left[\rho(\Aa_s)+(n-1)\Wa_s\right]ds\right)}:=\overline{\Va}_t/\Va_t .
$$

\end{theo}

The proof of this theorem is rather technical thus it is housed in Section~\ref{proof-gronwall} in the appendix.

\begin{cor}

When $\Ua_t=0=\Va_t$ we have
\begin{equation}\label{hilbert}
\displaystyle\EE\left(\Vert\Xa_t\Vert^{n}~\vert~\Fa_0\right)\leq 
\EE\left(\exp{\left(\int_0^t\left[n~\rho(\Aa_s)~+\frac{n(n-1)}{2}~ \Wa_s\right]~ds\right)}~|~\Fa_0\right)^{1/2}~\Vert\Xa_0\Vert^{n}.
\end{equation}

When $\rho(\Aa_t)\leq -a_t$ and $\Wa_t\leq w_t$ for some constants $a_t,w_t$, and $\Xa_0=0$ we have
\begin{equation}\label{hilbert-G}
\begin{array}{l}
\EE\left(\Vert\Xa_t\Vert^{n}\right)^{2/n}
\displaystyle\leq \displaystyle
\int_0^t~\exp{\left(-\left[
\int_s^t\lambda_n(a_u,w_u)~du+\frac{n-1}{2}\int_0^sw_u~du\right]\right)}~\\
\\
\hskip7cm\times \displaystyle\left[
 \EE\left(\Ua_s^n\right)^{1/n}+\frac{(n-1)}{2}~\EE\left(\Va_s^n\right)^{1/n}\right]~ds
\end{array}\end{equation}
with
$$
\lambda_n(a_s,w_s):=a_s-\frac{n-1}{2}~w_s.
$$
\end{cor}

\proof
The first assertion is a direct consequence of the estimates stated in Theorem~\ref{gronwall-2dim}.
Replacing  $\Wa_t$ and $\rho(\Aa_t)$ by $w_t$ and $(-a_t)$ from the start in the proof of Theorem~\ref{gronwall-2dim}
we find that

$$
\begin{array}{l}
\EE\left(\Vert\Xa_t\Vert^{n}~|~\Fa_0\right)^{2/n}
\displaystyle\leq \exp{\left(-\int_0^t\lambda_n(a_s,w_s)~ds\right)}~\Vert\Xa_0\Vert^{2}\\
\\
\hskip3cm+\displaystyle
\int_0^t~\exp{\left(-\left[
\int_s^t\lambda_n(a_u,w_u)~du+\frac{n-1}{2}\int_0^sw_u~du\right]\right)}~\\
\\
\hskip7cm\times \displaystyle\left[
 \EE\left(\Ua_s^n~|~\Fa_0\right)^{1/n}+\frac{(n-1)}{2}~\EE\left(\Va_s^n~|~\Fa_0\right)^{1/n}\right]~ds.
\end{array}$$
In the above display we have used the fact that
$$
\begin{array}{l}
\displaystyle\int_s^t\left(-a_u+\frac{n-1}{2}~w_u\right)~du+\int_0^s\left(-a_u+\frac{n-1}{2}~w_u\right)~du+
\int_0^s\left(a_u-(n-1)~w_u\right)~du\\
\\
\displaystyle\leq 
-\int_s^t\left(a_u-\frac{n-1}{2}~w_u\right)~du-\frac{n-1}{2}\int_0^sw_u~du.
\end{array}
$$
This ends the proof of the corollary.

\cqfd

\begin{prop}\label{prop-Laplace-Chi2}
Assume that $\rho(\Aa_t)\leq -a$ for some parameter $a>0$, and $\Xa_0=0=\Wa_t$. 
Also assume that for any $n\geq 1$ and any $t\geq 0$ we have
$$
 \EE\left(\Ua_t^n\right)^{1/n}\leq u_t\quad\mbox{and}\quad\EE\left(\Va_t^n\right)^{1/n}\leq v_t
$$
for some functions $u_t,v_t\geq 0$. In this situation, for any $\epsilon\in ]0,1]$ we have the uniform estimates
\begin{equation}\label{expo-chi2-unif}
\sup_{t\geq 0}{\EE\left(\exp{\left[\frac{(1-\epsilon)}{e}~\frac{1}{2v_t(a)}~\Vert\Xa_t\Vert^{2}\right]}\right)}\leq \frac{1}{2}~
\exp{\left(\frac{1-\epsilon}{e} \frac{u_t(a)}{v_t(a)}\right)}+\frac{e}{2\sqrt{2}}~\frac{1}{\sqrt{\epsilon}}
\end{equation}
for any functions $(u_t(a),v_t(a))$ such that
$$
\int_0^t~e^{-a(t-s)}~u_s~ds\leq u_t(a)\quad\mbox{and}\quad
\int_0^t~e^{-a(t-s)}~v_s~ds\leq v_t(a) .
$$
In addition, when $v_t=v$
for any $\epsilon\in[0,1]$ we have
\begin{eqnarray}\label{expo-chi2}
\EE\left(\displaystyle\exp{\left[\frac{a^2}{4v}~\epsilon~\int_0^t~\Vert\Xa_s\Vert^2~ds\right]}\right)
&\leq&\EE\left(\exp{\left[\frac{a}{v}~\frac{\epsilon}{1+\sqrt{1-\epsilon}}~
\int_0^t\Ua_s~ds\right]}\right)^{1/2} .
\end{eqnarray}
\end{prop}

The proof of the proposition is provided in the appendix, Section~\ref{proof-prop-appendix}.

We end this section with some comments on the estimate (\ref{expo-chi2}). Let us suppose that
$$
d\Vert\Xa_t\Vert^2=\left[-a~\Vert\Xa_t\Vert^2+u\right]~dt+d\Ma_t
$$
for some $u\geq 0$ (with $\Xa_0=0$).
In this case, by Jensen's inequality we have
$$
\begin{array}{l}
a~\EE(\Vert\Xa_t\Vert^2)=u~\left(1-e^{-at}\right)\\
\\
\displaystyle\Rightarrow
\EE\left(\displaystyle\exp{\left[\frac{a^2}{v}~\epsilon~\int_0^t~\Vert\Xa_s\Vert^2~ds\right]}\right)\geq 
\exp{\left[\frac{a}{v}~\epsilon~u~t~\left(1-\frac{1}{at}\left[1-e^{-at}\right]\right)\right]}\geq \exp {\left[\frac{a}{v}~\epsilon~u~t~
\left(1-\frac{1}{at}\right)\right]} .
\end{array}
$$
The r.h.s. of (\ref{expo-chi2}) gives the estimate
$$
\EE\left(\displaystyle\exp{\left[\frac{a^2}{v}~\epsilon~\int_0^t~\Vert\Xa_s\Vert^2~ds\right]}\right)\leq 
\exp{\left[\frac{a}{v}~\epsilon~u~t~\frac{1}{1+\sqrt{1-\epsilon}}\right]} .$$
The above estimates coincides for any $\epsilon\in[0,1]$ and any $u\geq 0$ as soon as
$$
t\geq \frac{1}{a}~\left(1+\frac{1}{\sqrt{1-\epsilon}}\right) .
$$

\section{Stochastic stability properties}\label{stab-prop-section}
\subsection{The signal process}\label{signal-sec}
This section is mainly concerned with the stochastic stability properties of the signal process.
One natural way to derive some useful concentration inequalities is to compare the flow of the 
stochastic process with the one of the noise free deterministic system discussed in the beginning of
Section~\ref{statement-section}.

We start with a brief review on the long-time behaviour of the semigroup $\varphi_t(x)$. It is readily checked that
\begin{eqnarray*}
\partial_t\Vert \varphi_t(x)-\varphi_t(y)\Vert^2
&\leq &-2\lambda_A\Vert \varphi_t(x)-\varphi_t(y)\Vert^2\Rightarrow
\Vert \varphi_t(x)-\varphi_t(y)\Vert\leq e^{-\lambda_At}~\Vert x-y\Vert .
\end{eqnarray*}
This contraction property ensures the existence and the uniqueness of a fixed point
$$
\forall t\geq 0\quad \varphi_t(x_{\star}):=x_{\star}\Longleftrightarrow
A(x_{\star})=0\Longrightarrow \Vert \varphi_t(x)-x_{\star}\Vert\leq e^{-\lambda_At}~\Vert x-x_{\star}\Vert . 
$$
We let $\delta\phi_t(x)$ be the Jacobian of the stochastic flow $\phi_t(x)$. We have
the matrix valued equation
$$
\partial_t \delta\phi_t(x) =\partial A(\phi_t(x))~\delta\phi_t(x)\Rightarrow
\delta\phi_t(x)~u=\exp{\left(\int_0^t\partial A(\phi_s(x))~ds\right)}~u
$$
for any $u\in\RR^{r_1}$. This implies that
$$
\Vert \delta\phi_t(x)\Vert:=\sup_{\Vert u\Vert\leq 1}
\Vert \delta\phi_t(x)~u\Vert\leq \exp{\left(-\lambda_{\partial A}t/2\right)}~\longrightarrow_{t\rightarrow\infty}~0 .
$$
Using the formula
$$
\phi_t(y)-\phi_t(x)=\int_0^1~\delta\phi_t(x+\epsilon(y-x))~(y-x)~d\epsilon ,
$$
we easily check the almost sure exponential stability property
\begin{equation}\label{stab-flow}
\Vert \phi_t(x)-\phi_t(y)\Vert\leq \exp{\left(-\lambda_{\partial A}t/2\right)}~\Vert x-y\Vert .
\end{equation}
The same analysis applies to estimate the Jacobian $\delta\varphi_t(x)$ of the deterministic flow $\varphi_t(x)$.
Using the estimate
$$
\Vert\phi_t(X_0)-\phi_t\left(\EE\left(X_0\right)\right)\Vert\leq e^{-\lambda_{\partial A}t/2}~
\Vert X_0-\EE\left(X_0\right)\Vert
$$
we also have
$$
\lambda_{\partial A}
\int_0^t~\Vert\phi_s(X_0)-\phi_s\left(\EE\left(X_0\right)\right)\Vert^2~ds\leq
\Vert X_0-\EE\left(X_0\right)\Vert^2
$$
from which we conclude that
$$
(\ref{chi-2})\Longrightarrow
\EE\left(\exp{\left(
\frac{\lambda_{\partial A}}{\chi(P_0)}~\int_0^t~\Vert\phi_s(X_0)-\phi_s\left(\EE\left(X_0\right)\right)\Vert^2~ds
\right)}\right)\leq e .
$$
The next proposition quantifies the relative stochastic stability of the flows $(\varphi_t,\phi_t)$
in terms of $\LL_n$-norms and $\chi$-square uniform Laplace estimates. 
\begin{prop}\label{prop-laplace-signal-dyns}
For any $n\geq 1$ and any $x\in\RR^{r_1}$ we have the uniform moment estimates
\begin{equation}\label{m-X-x}
\EE\left(\Vert \phi_t(x)-\varphi_t(x)\Vert^{2n}\right)^{1/n}\leq (n-1/2)~\tr(R_1)/\lambda_A .
\end{equation}
In addition, for any $\epsilon\in ]0,1]$ we have the uniform Laplace estimates
$$
\sup_{t\geq 0}{\EE\left(\exp{\left[\frac{(1-\epsilon)}{4e}~\frac{\lambda_A}{\tr(R_1)}~\Vert  \phi_t(x)-\varphi_t(x)\Vert^{2}\right]}\right)}\leq\frac{e}{2\sqrt{2}}
\frac{1}{\sqrt{\epsilon}}+\frac{1}{2}
\exp{\left[\frac{1-\epsilon}{4e}\right]}
$$
as well as
$$
\EE\left(\displaystyle\exp{\left[\frac{\lambda_A^2}{4\tr(R_1)}~\epsilon~\int_0^t~\Vert  \phi_s(x)-\varphi_s(x)\Vert^{2}~ds\right]}\right)
\leq\exp{\left[\frac{\lambda_A}{2}~\epsilon~t\right]} .
$$

\end{prop}

Combining (\ref{m-X-x}) with the concentration inequality (\ref{event-control})
we prove that  the probability of the
event
$$
\Vert \phi_t(x)-\varphi_t(x)\Vert^2\leq  \varpi(\delta)~\tr(R_1)/\lambda_A
$$
is greater than $1-e^{-\delta}$, for any $\delta\geq 0$ and any initial states 
$x\in\RR^{r_1}$.
This ends the proof of~(\ref{concentration-stab-signal}).

{\bf Proof of Proposition~\ref{prop-laplace-signal-dyns}:}

We have
$$
d(X_t-x_t)=[A(X_t)-A(x_t)]~dt+R^{1/2}_1dW_t
$$
with $X_0=x_0$,
and therefore
\begin{eqnarray*}
d\Vert X_t-x_t\Vert^2&=&\left[2\langle A(X_t)-A(x_t),X_t-x_t,\rangle+\tr(R_1)\right]~dt+dM_t\\
&\leq& \left[-2\lambda_A\Vert X_t-x_t\Vert^2+\tr(R_1)\right]dt+dM_t
\end{eqnarray*}
with the martingale
$$
dM_t:=2
\langle X_t-x_t,R^{1/2}_1dW_t\rangle
\Longrightarrow \partial_t\langle M\rangle_t=4~\mbox{\rm tr}\left(R_1(X_t-x_t)(X_t-x_t)^{\prime}\right)\leq 4~\tr(R_1)~\Vert X_t-x_t\Vert^2 .
$$
The end of the proof is now a direct consequence of (\ref{hilbert-G}) and 
 Proposition~\ref{prop-Laplace-Chi2} applied to
$$
\Xa_t=\Vert X_t-x_t\Vert , \qquad
\Aa_t x:=-ax=-2\lambda_A x ,\qquad \Ua_t=u=\tr(R_1)\quad\mbox{\rm and}\quad
\Va_t=v=4~\tr(R_1).
$$
The proof of the proposition is now completed.
\cqfd

\subsection{The Extended Kalman-Bucy filter}\label{lp-section}

This section is mainly concerned with the stochastic stability and the concentration properties of the semigroup of the EKF stochastic process. As for the signal process discussed in Section~\ref{signal-sec} 
these properties are related to $\LL_p$-mean error estimates and related $\chi$-square type Laplace inequalities.
Our main results are described by the following theorem.

Let $(\widehat{X}_t,P_t)$ be the solution of the evolution equations \eqref{EKF} starting at $(\widehat{X}_0,P_0)$.
\begin{theo}\label{theo-XwX}
For any $n\geq 1$ we have
\begin{equation}\label{moments-X-wX}
\displaystyle 
 \EE\left(\Vert\phi_t(\widehat{X}_0)-\widehat{X}_t\Vert^{n}\right)^{2/n}
\displaystyle \leq (2n-1)\left\{\frac{\tr(R_1)}{\lambda_A}~\frac{\sigma^2_{\partial A}}{2}+
\frac{\vert e^{-\lambda_At}- e^{-\lambda_{\partial A}t}\vert}{\vert \lambda_A-\lambda_{\partial A}\vert}~
\rho(S)~\tr(P_0)^2 \right\} .
\end{equation}
For any $\epsilon\in ]0,1]$ and any $P_0$ there exists some time horizon $t_0(\epsilon,P_0)$ such that
\begin{equation}\label{expo-chi2-unif-XwX}
\displaystyle \sup_{t\geq t_0(\epsilon,P_0)}{\EE\left(\exp{
\left[
\frac{(1-\epsilon)}{4e\sigma^2_{\partial A}}~\frac{\lambda_A}{\tr(R_1)}~
~\Vert\phi_t(\widehat{X}_0)-\widehat{X}_t\Vert^{2}
 \right]}
\right)}
\leq
\frac{1}{2}~
\exp{\left(\frac{1-\epsilon}{4e}\right)}+\frac{e}{2\sqrt{2}}~\frac{1}{\sqrt{\epsilon}} .
\end{equation}
In addition for any $t\geq s\geq 0$ and any $\epsilon\in ]0,1]$ we have
\begin{eqnarray}\label{expo-chi2-XwX}
\EE\left(
\displaystyle\exp{\left[\frac{\epsilon}{1+\pi_{\partial A}(s)}~\frac{\lambda_A^2}{4\tr(R_1)}~
\int_s^t~\Vert\phi_{r-s}(\widehat{X}_s)-\widehat{X}_r\Vert^2~ds\right]}\right)
&\leq&\exp{\left[\frac{\lambda_A}{2}~\epsilon~(t-s)
\right]}. 
\end{eqnarray}

\end{theo}

Before getting into the details of the proof of this theorem we mention that (\ref{concentration-stab})
is a direct consequence of 
 (\ref{moments-X-wX}) combined
 with (\ref{event-control}) and (\ref{stab-flow}). Indeed, applying (\ref{event-control}) to
 $$
 Z=\Vert\phi_t(\widehat{X}_0)-\widehat{X}_t\Vert
 \quad\mbox{\rm and}
 \quad
 z^2=4\left\{\frac{\tr(R_1)}{\lambda_A}~\frac{\sigma^2_{\partial A}}{2}+
\frac{\vert e^{-\lambda_At}- e^{-\lambda_{\partial A}t}\vert}{\vert \lambda_A-\lambda_{\partial A}\vert}~
\rho(S)~\tr(P_0)^2 \right\}
 $$
 by (\ref{moments-X-wX}) we readily check that 
 the probability of the
events
$$
\begin{array}{l}
\Vert \phi_t(x)-\Phi_t(\widehat{x},p)\Vert^2\\
\\
\displaystyle\leq 2
\exp{\left(-\lambda_{\partial A}t\right)}~\Vert x-\widehat{x}\Vert^2+8\varpi(\delta)~\left\{\frac{\tr(R_1)}{\lambda_A}~\frac{\sigma^2_{\partial A}}{2}+
\frac{\vert e^{-\lambda_At}- e^{-\lambda_{\partial A}t}\vert}{\vert \lambda_A-\lambda_{\partial A}\vert}~
\rho(S)~\tr(P_0)^2 \right\}
\end{array}
$$
is greater than $1-e^{-\delta}$, for any $\delta\geq 0$ and any initial states 
$(x,\widehat{x},p)\in\RR^{r_1+r_1+(r_1\times r_1)}$. In this connection, the Laplace estimates (\ref{expo-chi2-XwX})
readily imply that  the probability of the
events
$$
\frac{1}{t-s}
\int_s^t~\Vert\phi_u(\widehat{X}_s)-\widehat{X}_u\Vert^2~du\leq \left(\frac{1}{2}+\frac{\delta}{\lambda_A}\right)~
(1+\pi_{\partial A}(s))~\frac{4\tr(R_1)}{\lambda_A}
$$
 is greater than $1-e^{-\delta}$, for any $\delta\geq 0$ and any time horizon $t$.

Now we come to the proof of the theorem.

 {\bf Proof of Theorem~\ref{theo-XwX}:}\\
We set
$
\Xa_t:=\phi_t(\widehat{X}_0)-\widehat{X}_t
$. We have
$$
d\Vert\Xa_t\Vert^2\leq \left(2~\langle A(\phi_t(\widehat{X}_0))-A(\widehat{X}_t),\Xa_t\rangle-2~\langle P_tS\Xa_t,\Xa_t\rangle+\tr(R_1)+\tr(SP^2_t)\right)
~dt+d\Ma_t
$$
with the martingale
$$
d\Ma_t:=2\Xa_t^{\prime}\left(
R^{-1/2}_{1}dW_t
-P_tB^{\prime}R^{-1/2}_{2}dV_t\right)\rightarrow\partial_t\langle \Ma_t\rangle_t .
$$
This yields the estimate
$$
d\Vert\Xa_t\Vert^2\leq \left[-2\lambda_A~\Vert\Xa_t\Vert^2+\Ua_t\right]
~dt+d\Ma_t
\quad\mbox{\rm with}
\quad
\Ua_t=u_t:=\left[\tr(R_1)+ \tau_t^2(P)~\rho(S)\right] .
$$
Also observe that
\begin{eqnarray*}
\partial_t\langle \Ma\rangle_t&\leq& 4~\Vert\Xa_t\Vert^2~\left(\tr(R_1)+\tr(SP_t^2)\right)\leq
\Va_t~\Vert\Xa_t\Vert^2
\quad\mbox{\rm
with $\Va_t=v_t:=4u_t$.
}
\end{eqnarray*}
On the other hand we have
$$
2e^{-2\lambda_At}\int_0^t~e^{2(\lambda_A-\lambda_{\partial A})s}~ds=\frac{\vert e^{-\lambda_At}- e^{-\lambda_{\partial A}t}\vert}{\vert \lambda_A-\lambda_{\partial A}\vert}~.
$$
This implies that
\begin{eqnarray*}
\int_0^t~e^{-2\lambda_A(t-s)}~\tau_s^2~ds&\leq&~2\tr(P_0)^2 ~e^{-2\lambda_At}
\int_0^t~e^{-2\Delta_As}~ds+\frac{1}{\lambda_A}~\left(\frac{\tr(R_1)}{\lambda_{\partial A}}\right)^2\\
&\leq&\frac{\vert e^{-\lambda_At}- e^{-\lambda_{\partial A}t}\vert}{\vert \lambda_A-\lambda_{\partial A}\vert}~
~\tr(P_0)^2~+~\frac{1}{\lambda_A}~\left(\frac{\tr(R_1)}{\lambda_{\partial A}}\right)^2 .
\end{eqnarray*}
This implies that
\begin{eqnarray*}
\int_0^t~e^{-2\lambda_A(t-s)}~u_s~ds&\leq& \frac{\tr(R_1)}{2\lambda_A}+\rho(S)\left[
\frac{\vert e^{-\lambda_At}- e^{-\lambda_{\partial A}t}\vert}{\vert \lambda_A-\lambda_{\partial A}\vert}~
~\tr(P_0)^2~+~\frac{1}{\lambda_A}~\left(\frac{\tr(R_1)}{\lambda_{\partial A}}\right)^2
\right]\\
&=& \frac{\tr(R_1)}{\lambda_A}\left(\frac{1}{2}+\frac{\rho(S)}{\lambda_{\partial A}}~
\frac{\tr(R_1)}{\lambda_{\partial A}}\right)+
\frac{\vert e^{-\lambda_At}- e^{-\lambda_{\partial A}t}\vert}{\vert \lambda_A-\lambda_{\partial A}\vert}~
\rho(S)~\tr(P_0)^2~\\
&:=&u_t(a):=v_t(a)/4 .
\end{eqnarray*}

Applying 
 Proposition~\ref{prop-Laplace-Chi2} to 
 $
\Aa_t x:=-ax=-2\lambda_A x
$,
we find that
$$
{\EE\left(\exp{\left[\frac{(1-\epsilon)}{e}~\frac{1}{8u_t(a)}~\Vert\Xa_t\Vert^{2}\right]}\right)}\leq \frac{1}{2}~
\exp{\left(\frac{1-\epsilon}{4e}\right)}+\frac{e}{2\sqrt{2}}~\frac{1}{\sqrt{\epsilon}}
\quad\mbox{\rm
for any $\epsilon\in ]0,1]$.}
$$

Using the fact that for any non negative real numbers $x,y,\lambda$ we have
\begin{eqnarray*}
\frac{1}{x+e^{-\lambda}~y}&=&\frac{1}{x}\left(1-\frac{e^{-\lambda}~y/x}{1+e^{-\lambda}~y/x}\right)\geq
\frac{1}{x}\left(1-e^{-\lambda}~y/x\right) 
\end{eqnarray*}
and
$$
\frac{\vert e^{-\lambda_At}- e^{-\lambda_{\partial A}t}\vert}{\vert \lambda_A-\lambda_{\partial A}\vert}~\longrightarrow_{t\rightarrow\infty}~0 ,
$$
we find that
\begin{eqnarray*}
\frac{1}{u_t(a)}&\geq&  (1-\epsilon)~ \frac{\lambda_A}{\lambda_{\partial A}}~\left[
\frac{\tr(R_1)}{\lambda_{\partial A}}\left(\frac{1}{2}+\frac{\rho(S)}{\lambda_{\partial A}}~
\frac{\tr(R_1)}{\lambda_{\partial A}}\right)\right]^{-1}
\end{eqnarray*}
for any $t\geq t(\epsilon)$, for any $\epsilon\in [0,1[$ and some $t(\epsilon)$.

 The end of the proof is now a direct consequence of (\ref{hilbert-G}) and 
 Proposition~\ref{prop-Laplace-Chi2} applied to 
 $$
\Aa_t x:=-ax=-2\lambda_A x\quad \mbox{\rm and}\quad u_t=v_t/4\leq 
u=v/4:=\left[\tr(R_1)+ \tau_s^2(P)~\rho(S)\right]
$$
with $t\in [s,\infty[$.
The proof of the theorem is now completed.
\cqfd

\section{Proof of Theorem~\ref{stab-theo}}

We let $(\widehat{X}_t,P_t)$ be the solution of Equations (\ref{EKF}) starting at $(\widehat{X}_0,P_0)$. We denote by  $(\widecheck{X}_t,\widecheck{P}_t)$ the solution of these equations starting at 
some possibly different state $(\widecheck{X}_0,\widecheck{P}_0)$. 
Firstly we have
$$
\left((\ref{riccati-bound}),~
(\ref{m-X-x})~\mbox{\rm and}~
(\ref{moments-X-wX})\right)~\Longrightarrow\forall n\geq 1\qquad~\sup_{t\geq 0}{\EE\left( [\Vert \widehat{X}_t-\widecheck{X}_t
\Vert^2+\Vert P_t-\widecheck{P}_t\Vert_F^2]^n   
 \right)}<\infty .
$$

We couple the equations with the same observation processes. In this situation we find the 
evolution equation
\begin{eqnarray*}
d(\widehat{X}_t-\widecheck{X}_t)&=&\left(A(\widehat{X}_t)-A(\widecheck{X}_t)\right)dt+P_tB^{\prime}R^{-1}_{2}~\left[dY_t-B\widehat{X}_tdt\right]\\
&&\hskip6cm-\widecheck{P}_tB^{\prime}R^{-1}_{2}~\left[dY_t-B\widecheck{X}_tdt\right]
\\
&&\\
&=&\left(\left[A(\widehat{X}_t)-A(\widecheck{X}_t)\right]+\widecheck{P}_t B^{\prime}R^{-1}_{2}~\left[B~(\widecheck{X}_t-X_t)\right]\right)dt\\
\\
&&\hskip.3cm
+(P_t-\widecheck{P}_t)B^{\prime}R^{-1}_{2}~\left[B(X_t-\widehat{X}_t)\right]~dt+\widecheck{P}_tB^{\prime}R^{-1}_{2}~\left[B(X_t-\widehat{X}_t)\right]dt+dM_t 
\end{eqnarray*}
with the martingale
$$
dM_t:=\left[P_t-\widecheck{P}_t\right]B^{\prime}R^{-1/2}_{2}~dV_t .
$$
This implies that
$$
\begin{array}{l}
d(\widehat{X}_t-\widecheck{X}_t)\\
\\
=\left([A(\widehat{X}_t)-A(\widecheck{X}_t)]+
\widecheck{P}_t S~(\widecheck{X}_t-\widehat{X}_t)
+(P_t-\widecheck{P}_t)S(X_t-\widehat{X}_t)\right)~dt+dM_t
\end{array}
$$
from which we conclude that
\begin{equation}\label{eq-diff-norm}
\begin{array}{l}
d\Vert\widehat{X}_t-\widecheck{X}_t\Vert^2\\
\\
=2\langle \widehat{X}_t-\widecheck{X}_t,[A(\widehat{X}_t)-A(\widecheck{X}_t)]-
\widecheck{P}_t S~(\widehat{X}_t-\widecheck{X}_t)+(P_t-\widecheck{P}_t)S(X_t-\widehat{X}_t)~\rangle~dt\\
\\
\hskip2cm
+\tr\left(S\left[P_t-\widecheck{P}_t\right]^2\right)
dt+d\overline{M}_t
\end{array}
\end{equation}
with the martingale
$$
d\overline{M}_t=2~\langle \widehat{X}_t-\widecheck{X}_t, dM_t\rangle .
$$
The angle bracket of $\overline{M}_t$ is given by
\begin{eqnarray*}
\partial_t\langle \overline{M}\rangle_t&=&4~\tr\left(\left[P_t-\widecheck{P}_t\right]S\left[P_t-\widecheck{P}_t\right](\widehat{X}_t-\widecheck{X}_t)(\widehat{X}_t-\widecheck{X}_t)^{\prime}\right)\\
&=&4~\langle (\widehat{X}_t-\widecheck{X}_t),\left[P_t-\widecheck{P}_t\right]S\left[P_t-\widecheck{P}_t\right](\widehat{X}_t-\widecheck{X}_t)\rangle\\
&\leq& 4~\rho(S)~\Vert\widehat{X}_t-\widecheck{X}_t\Vert^2~\Vert P_t-\widecheck{P}_t\Vert^2_F\leq 
2~\rho(S)~ \left(\Vert\widehat{X}_t-\widecheck{X}_t\Vert^2+\Vert P_t-\widecheck{P}_t\Vert^2_F\right)^2 .
\end{eqnarray*}
Recalling that $\lambda_A\geq \lambda_{\partial A}/2$, also observe that the drift term in (\ref{eq-diff-norm}) is bounded by
$$
\begin{array}{l}
-\lambda_{\partial A}~\Vert\widehat{X}_t-\widecheck{X}_t\Vert^2+2\beta_t~\Vert\widehat{X}_t-\widecheck{X}_t\Vert
\Vert P_t-\widecheck{P}_t\Vert_F
+\rho(S)~\Vert P_t-\widecheck{P}_t\Vert^2_F
\end{array}$$
with 
$$
\beta_t:=\rho(S)\Vert X_t-\widehat{X}_t\Vert .
$$

In much the same way we have
$$
\begin{array}{l}
\partial_t(P_t-\widecheck{P}_t)\\
\\=\left(\partial A(\widehat{X}_t)P_t-\partial A(\widecheck{X}_t)\widecheck{P}_t\right)+\left(\partial A(\widehat{X}_t)P_t-\partial A(\widecheck{X}_t)\widecheck{P}_t\right)^{\prime}+\widecheck{P}_t~S\widecheck{P}_t
-P_t~SP_t\\
\\
=\left([\partial A(\widehat{X}_t)-\partial A(\widecheck{X}_t)]P_t+\partial A(\widecheck{X}_t)
[P_t-\widecheck{P}_t]\right)+\left([\partial A(\widehat{X}_t)-\partial A(\widecheck{X}_t)]P_t+\partial A(\widecheck{X}_t)
[P_t-\widecheck{P}_t]\right)^{\prime}\\
\\
\hskip2cm+\frac{1}{2}~(\widecheck{P}_t+P_t)S(\widecheck{P}_t-P_t)+\frac{1}{2}~(\widecheck{P}_t-P_t)S(\widecheck{P}_t+P_t) .
\end{array}
$$
In the last assertion we have used the matrix decomposition
$$
PSP-QSQ=\frac{1}{2}~(P+Q)S(P-Q)+\frac{1}{2}~(P-Q)S(P+Q) .
$$
Recalling that
$$
2^{-1}\partial_t\Vert P_t-\widecheck{P}_t\Vert^2
=2^{-1}\partial_t\langle P_t-\widecheck{P}_t,P_t-\widecheck{P}_t\rangle
=\langle P_t-\widecheck{P}_t,\partial_t(P_t-\widecheck{P}_t)\rangle=\tr\left((P_t-\widecheck{P}_t)\partial_t(P_t-\widecheck{P}_t)\right) ,
$$
we find that
\begin{eqnarray*}
2^{-1}\partial_t\Vert P_t-\widecheck{P}_t\Vert_F^2&=&
2\tr\left(\partial A(\widecheck{X}_t)(P_t-\widecheck{P}_t)^2\right)+2\tr\left([\partial A(\widehat{X}_t)-\partial A(\widecheck{X}_t)]P_t(P_t-\widecheck{P}_t)\right)\\
&&\hskip6cm-\tr\left(
(\widecheck{P}_t+P_t)S(\widecheck{P}_t-P_t)^2
\right)\\
&\leq &2\tr\left(\partial A(\widecheck{X}_t)(P_t-\widecheck{P}_t)^2\right)+2\tr\left([\partial A(\widehat{X}_t)-\partial A(\widecheck{X}_t)]P_t(P_t-\widecheck{P}_t)\right) .
\end{eqnarray*}
This implies that
\begin{eqnarray*}
\partial_t\Vert P_t-\widecheck{P}_t\Vert_F^2
&\leq& -2\lambda_{\partial A}~\Vert P_t-\widecheck{P}_t\Vert_F^2+2\alpha_{t}\Vert P_t-\widecheck{P}_t\Vert_F\Vert \widehat{X}_t-\widecheck{X}_t\Vert
\end{eqnarray*}
with
$$
\alpha_t:=2\kappa_{\partial A}\tau_t( P) .
$$

We set
$$
\Xa_t=\left(
\begin{array}{c}
\Vert \widehat{X}_t-\widecheck{X}_t\Vert\\
\Vert P_t-\widecheck{P}_t\Vert_F
\end{array}
\right)\in \Ha:=\RR^2
\Longrightarrow
d\Vert\Xa_t\Vert^2\leq \langle\Xa_t,\Aa_t\Xa_t\rangle~dt+d\Ma_t
$$
with
$$
\Aa_t=\left(
\begin{array}{cc}
-\lambda_{\partial A}&2\beta_t\\
2\alpha_{t}&-2\lambda_{\partial A}+\rho(S)
\end{array}
\right)\quad\mbox{\rm and}\quad
\Ma_t=
\overline{M}_t .
$$
Notice that
$$
(\Aa_t+\Aa_t^{\prime})/2=\left(
\begin{array}{cc}
-\lambda_{\partial A}&\beta_t+\alpha_t\\
\beta_t+\alpha_t&-2\lambda_{\partial A}+\rho(S)
\end{array}
\right) .
$$

Observe that
\begin{eqnarray*}
\rho(\Aa_t)&:=&
\lambda_{\tiny max}\left((\Aa_t+\Aa_t^{\prime})/2\right)\\
&=&-\frac{1}{2}\left(3\lambda_{\partial A}-\rho(S)\right)+\sqrt{\frac{1}{4}
\left(\lambda_{\partial A}-\rho(S)\right)^2+\left(\beta_t+\alpha_t\right)^2
}\leq -\lambda_{\partial A}+\beta_t+\alpha_t\\
&\leq& \overline{\rho}(\Aa_t):=-\left(\lambda_{\partial A}-2\kappa_{\partial A}\tau_t(P)-\rho(S)~\Vert X_t-\widehat{X}_t\Vert\right)
\end{eqnarray*}

The final step is based on the following technical lemma.

\begin{lem}\label{techn-lemma}
Assume Condition  (\ref{regularity-conditon-HAS-0})  is satisfied for some  $\alpha> 1$. 
In this situation, for any $\epsilon\in ]0,1]$  there exists some time horizon $s$ such that
for any $t\geq s$ we have the almost sure estimate
\begin{equation}\label{estimation-lemma}
\delta:=\frac{1}{2}~ 
\sqrt{\frac{\lambda_{\partial A}}{\rho(S)}}\Longrightarrow\begin{array}[t]{l}
\displaystyle\EE\left(
\exp{\left[
\delta
\int_s^t\left\{
(\overline{\rho}(\Aa_u)+(\delta-1)\rho(S))
\right\}~du
\right]
}~|~\Fa_s\right)^{1/\delta}\\
\\
\hskip3cm\leq \displaystyle \Za_{s}~\displaystyle \exp{\left(-\left(1-\epsilon\right)\Lambda(t-s)\right)}
\end{array}
\end{equation}
for some positive random process $\Za_{t}$ s.t. 
$
\sup_{t\geq 0}{\EE\left(\Za_t^{\alpha\delta}\right)}<\infty
$.

\end{lem}

The end of the proof of Theorem~\ref{stab-theo} is  a direct consequence of this lemma, so we give it first.
Combining (\ref{estimation-lemma}) with (\ref{gronwall-n+}) we find that
$$
\EE\left(\Vert\Xa_t\Vert^{\delta}~|~\Fa_s\right)^{2/\delta}
\displaystyle\leq  \Za_{s}~\displaystyle \exp{\left(-\left(1-\epsilon\right)\Lambda(t-s)\right)}~\Vert\Xa_s\Vert^2
$$
This ends the proof of Theorem~\ref{stab-theo}.

Now we come to the proof of the lemma.

{\bf Proof of Lemma~\ref{techn-lemma}:}

For any $t\geq s\geq 0$ we have the estimate
\begin{eqnarray*}
\overline{\rho}(\Aa_t)
&\leq&
-\left(\lambda_{\partial A}-2\kappa_{\partial A}\tau_s(P)\right)+\beta_t=-\Delta_{\partial A}(s)+\rho(S)~\Vert X_t-\widehat{X}_t\Vert
\end{eqnarray*}
with
$$
\Delta_{\partial A}(s):=
\lambda_{\partial A}-2\kappa_{\partial A}~\tau_s(P)\longrightarrow_{s\rightarrow\infty}
\Delta_{\partial A}:=
\lambda_{\partial A}-2\kappa_{\partial A}~\tr(R_1)/\lambda_{\partial A}>0
$$ 
as soon as
$$
\lambda_{\partial A}>\sqrt{2\kappa_{\partial A}~\tr(R_1)} .
$$

For any $\epsilon\in ]0,1]$, there exists some time horizon $ \varsigma_{\epsilon}(P_0)$  
such that 
\begin{eqnarray*}
t\geq s\geq \varsigma_{\epsilon}(P_0)~&\Longrightarrow&~
(1-\epsilon)\leq \Delta_{\partial A}(s)/\Delta_{\partial A}\leq 1\\&\Longrightarrow&
\overline{\rho}(\Aa_t)\leq -(1-\epsilon)\Delta_{\partial A}+\rho(S)~\Vert X_t-\widehat{X}_t\Vert
\end{eqnarray*}
On the other hand, the contraction inequality \eqref{stab-flow} implies that
\begin{eqnarray*}
\int_s^t~ \Vert X_r-\widehat{X}_r\Vert~dr&=&\int_s^t~ \Vert \phi_{r-s}(\phi_s(X_0))-\widehat{X}_r\Vert~dr\\
&\leq &\int_s^t~ \Vert \phi_{r-s}(X_s)-\phi_{r-s}(\widehat{X}_s)\Vert~dr
+\int_s^t~ \Vert \phi_{r-s}(\widehat{X}_s)-\widehat{X}_r\Vert~dr\\
&\leq &\left(\int_s^t~e^{-\lambda_{\partial A}r/ 2}~dr\right)~ \Vert X_s-\widehat{X}_s\Vert+\int_s^t~ \Vert \phi_{r-s}(\widehat{X}_s)-\widehat{X}_r\Vert~dr\\
&\leq & 2 \Vert X_s-\widehat{X}_s\Vert/\lambda_{\partial A}+\int_s^t~ \Vert \phi_{r-s}(\widehat{X}_s)-\widehat{X}_r\Vert~dr
\end{eqnarray*}
 The above inequality yields the almost sure estimate 
$$
\begin{array}{l}
\displaystyle\EE\left(
\exp{\left[
\delta
\left\{
\int_s^t(\overline{\rho}(\Aa_u)+(\delta-1)\rho(S))
~du \right\}
\right]
}~|~\Fa_s\right)\\
\\
\leq \displaystyle
\exp{\left[
-\delta
\left\{
(1-\epsilon)\Delta_{\partial A}+(1-\delta)\rho(S)
\right\}(t-s)
\right]
}\\
\\
\hskip5cm\times\displaystyle ~\Za_s~
\EE\left(
\exp{\left[
\delta
\left\{\rho(S)
\int_s^t~ \Vert \phi_{u-s}(\widehat{X}_s)-\widehat{X}_u\Vert~
~du \right\}
\right]
}~|~\Fa_s\right)
\end{array}
$$
with
$$
\Za_s:=
\exp{\left[
2\rho(S)
\Vert X_s-\widehat{X}_s\Vert/\lambda_{\partial A}
\right]
} .
$$
Using the estimate $x-1/4\leq x^2$, which is valid for any $x$, we have
$$
\int_s^t~((\Vert \phi_{u-s}(\widehat{X}_s)-\widehat{X}_u\Vert- 1/4)+ 1/4)~du\leq (t-s)/4+\int_s^t~\Vert \phi_{u-s}(\widehat{X}_s)-\widehat{X}_u\Vert^2du .
$$
We find that
$$
\begin{array}{l}
\displaystyle\EE\left(
\exp{\left[
\delta
\left\{
\int_s^t(\overline{\rho}(\Aa_u)+(\delta-1)\rho(S))
\right\}~du
\right]
}~|~\Fa_s\right)\\
\\
\leq \displaystyle
\exp{\left[
-\delta
\left\{
(1-\epsilon)\Delta_{\partial A}+(3/4-\delta)\rho(S)
\right\}(t-s)
\right]
}\\
\\
\hskip5cm\times\displaystyle~\Za_s~
\EE\left(
\exp{\left[
\delta
\rho(S)
\int_s^t~ \Vert \phi_{u-s}(\widehat{X}_s)-\widehat{X}_u\Vert^2~du
\right]
}~|~\Fa_s\right)
\end{array}
$$

By (\ref{regularity-conditon-HAS-0}) we can also choose $s$ sufficiently large so that
$$
\begin{array}{l}
\displaystyle
\delta=\frac{1}{2}\sqrt{\frac{\lambda_{\partial A}}{\rho(S)}}\Longrightarrow\delta\rho(S)4\tr(R_1)~(1+\pi_{\partial A}(s))\leq~\lambda_A\frac{\lambda_{\partial A}}{2} .
\end{array}
$$
In this situation, by (\ref{expo-chi2-XwX}), we have
$$
\begin{array}{l}
\displaystyle\delta\rho(S)\leq~\overbrace{\frac{\lambda_{\partial A}}{2\lambda_A}}^{\leq 1}~\frac{1}{1+\pi_{\partial A}(s)}~\frac{\lambda_A^2}{4\tr(R_1)}\\
\\
\displaystyle\Rightarrow \begin{array}[t]{rcl}
\displaystyle\EE\left(\exp{\left(\delta\rho(S)~\int_s^t~\Vert \phi_{r-s}(\widehat{X}_s)-\widehat{X}_r\Vert^2~dr\right)}~|~\Fa_s\right)&\leq& \displaystyle\exp{\left[\frac{\lambda_{\partial A}}{4}~(t-s)\right]} .\\
\end{array}\end{array}
$$
We conclude that
$$
\begin{array}{l}
\displaystyle\EE\left(
\exp{\left[
\delta
\left\{
\int_s^t(\overline{\rho}(\Aa_u)+(\delta-1)\rho(S))
\right\}~du
\right]
}~|~\Fa_s\right)^{1/\delta}\\
\\
\leq \displaystyle
\exp{\left[
-(1-\epsilon)~\lambda_{\partial A}
\left\{
1-2~\frac{\kappa_{\partial A}}{\lambda_{\partial A}}~\frac{\tr(R_1)}{\lambda_{\partial A}}-\frac{1}{4\delta}+(3/4-\delta)~\frac{\rho(S)}{\lambda_{\partial A}}
\right\}(t-s)
\right]
}~\Za_s\\
\\\leq \displaystyle
\exp{\left[
-(1-\epsilon)~
\Lambda(t-s)
\right]
}~\Za_s .
\end{array}
$$
The last assertion comes from the formula
$$
\begin{array}{l}
\displaystyle
\delta=\frac{1}{2}\sqrt{\frac{\lambda_{\partial A}}{\rho(S)}}\\
\\
\displaystyle\Longrightarrow
1+\frac{3}{4}~\frac{\rho(S)}{\lambda_{\partial A}}-2~\frac{\kappa_{\partial A}}{\lambda_{\partial A}}~\frac{\tr(R_1)}{\lambda_{\partial A}}-\frac{1}{4\delta}-\delta~\frac{\rho(S)}{\lambda_{\partial A}}\\
\\
\displaystyle=
1-2~\frac{\kappa_{\partial A}}{\lambda_{\partial A}}~\frac{\tr(R_1)}{\lambda_{\partial A}}+~
\sqrt{\frac{\rho(S)}{\lambda_{\partial A}}}\left[\frac{3}{4}
\sqrt{\frac{\rho(S)}{\lambda_{\partial A}}}-1\right]>0 .
\end{array}
$$

On the other hand we have
\begin{eqnarray*}
\Vert X_s-\widehat{X}_s\Vert&=&\Vert \phi_s(X_0)-\widehat{X}_s\Vert\leq \Vert \phi_s(X_0)-\phi_s(\widehat{X}_0)\Vert+\Vert \phi_s(\widehat{X}_0)-\widehat{X}_s\Vert\\
&\leq& e^{-\lambda_{\partial A}s/2}~\Vert X_0-\EE\left(X_0\right)\Vert+\Vert \phi_s(\widehat{X}_0)-\widehat{X}_s\Vert .
\end{eqnarray*}
This shows that
$$
\Za_s\leq \exp{\left(\frac{\rho(S)}{\lambda_{\partial A}}~\Vert \phi_s(\widehat{X}_0)-\widehat{X}_s\Vert
\right)}~\exp{\left(\frac{\rho(S)}{\lambda_{\partial A}}~e^{-\lambda_{\partial A}s/2}~\Vert X_0-\EE\left(X_0\right)\Vert\right)} .
$$

Under the assumption (\ref{regularity-conditon-HAS-0}) and using (\ref{expo-chi2-unif-XwX}) we have
$$
\begin{array}{l}
\displaystyle \alpha~\delta~
\frac{\rho(S)}{\lambda_{\partial A}}=\frac{\alpha}{2}~\sqrt{\frac{\rho(S)}{\lambda_{\partial A}}}
< \frac{1}{8e \sigma^2_{\partial A}}~\frac{\lambda_{\partial A}}{\tr(R_1)} 
\\
\\
\Longrightarrow \exists p> 1~:~
\sup_{s\geq 0}{\EE\left(\exp{\left(p\alpha\delta \Vert \phi_s(\widehat{X}_0)-\widehat{X}_s\Vert (\rho(S)/\lambda_{\partial A}) \right)}\right)}<\infty .
\end{array}
$$
We can also choose $s$ sufficiently large so that
$$
\begin{array}{l}
\displaystyle\frac{ \alpha}{\alpha-1}~\delta\rho(S)~ e^{-\lambda_{\partial A}s/2}<\lambda_{\partial A}/\chi(P_0)\\
\\
\displaystyle\Longrightarrow 
\EE\left(
\exp{
\left(\alpha\delta\frac{p}{p-1}~(\rho(S)/\lambda_{\partial A})~ e^{-\lambda_{\partial A}s/2}~
\Vert X_0-\EE\left(X_0\right)\Vert\right)
}\right)\leq e.
\end{array}
$$
This ends the proof of the lemma.\cqfd
\section{Appendix}
\subsection{Concentration properties and Laplace estimates}

This section is mainly concerned with the proof of (\ref{chi-2}).

The initial state $X_0$ of the signal is a Gaussian  random variable
with mean $\widehat{X}_0$  and some covariance matrix $P_0$. In this case $X_0-\widehat{X}_0\stackrel{\tiny law}{=}P^{1/2}_0W_1$ and
$$
\EE\left(\Vert X_0-\widehat{X}_0\Vert^{2n}\right)\leq \rho(P_0)^n~\EE\left(\Vert W_1\Vert^{2n}\right) .
$$
Recalling that $\Vert W_1\Vert^{2}$ is distributed according to the chi-squared distribution with $r_1$ 
degrees of freedom we have
$$
\forall\gamma<1/(2\rho(P_0))\qquad
\EE\left(e^{\gamma \Vert X_0-\widehat{X}_0\Vert^{2}}\right)\leq \EE\left(e^{\gamma \rho(P_0)\Vert W_1\Vert^{2}}\right)=
\left(1-2\gamma\rho(P_0)\right)^{-r_1/2}<\infty .
$$
Using the fact that
$$
-t-\frac{1}{2}\log{(1-2t)}=t^2~\sum_{n\geq 0}~\frac{2}{2+n}~(2t)^n\leq \frac{t^2}{1-2t}\Rightarrow
\left(1-2t\right)^{-1/2}\leq \exp{\left(t+\frac{t^2}{1-2t}\right)}
$$
for any $0<t<1/2$, we check that
$$
\forall 0<t<(1-r_1\epsilon)/2\qquad
t+\frac{t^2}{1-2t}=t~\frac{1-t}{1-2t}\leq \frac{t}{r_1\epsilon}
$$
for any $\epsilon\in ]0,1/r_1[$. This yields
$$
\forall 0<\gamma<(1-r_1\epsilon)/(2\rho(P_0))\qquad
\EE\left(e^{\gamma \Vert X_0-\widehat{X}_0\Vert^{2}}\right)\leq\exp{
\left(\rho(P_0)\gamma/\epsilon
\right)} .
$$
Choosing
$
\gamma=(1-2r_1\epsilon)/(2\rho(P_0))
$, with $\epsilon\in ]0,1/(2r_1)[$ we find that
$$
\forall \epsilon\in ]0,1/(2r_1)[\qquad
\EE\left(\exp{\left[\left(\frac{1}{2}-(r_1\epsilon)\right) \frac{\Vert X_0-\widehat{X}_0\Vert^{2}}{\rho(P_0)}\right]}\right)\leq\exp{
\left(\left(\frac{1}{2}-r_1\epsilon\right)\frac{1}{\epsilon}
\right)} .
$$
We check (\ref{chi-2}) by choosing 
$$
\epsilon=\frac{2r_1-1}{4r_1^2}=\frac{1}{2r_1}-\frac{1}{4r_1^2}\Rightarrow
\EE\left(\exp{\left[ \frac{\Vert X_0-\widehat{X}_0\Vert^{2}}{4r_1\rho(P_0)}\right]}\right)\leq\exp{
\left(\frac{r_1}{2r_1-1}
\right)}\leq e .
$$
\cqfd

More generally, for any non negative random variable $Z$ such that
$$
\EE\left(Z^{2n}\right)^{1/(2n)}\leq z~\sqrt{n}\quad\mbox{\rm for some parameter $z\geq 0$}
$$
and for any $n\geq 1$ we have
$$
\EE\left(Z^{2n}\right)\leq (z^2n)^n\leq \frac{e}{\sqrt{2}}~\left(\frac{e}{2}~z^2\right)^{n} \EE(V^{2n})
$$
for some Gaussian and centered random variable $V$ with unit variance. We check this claim using Stirling approximation
\begin{eqnarray}
\EE(V^{2n})&=&2^{-n}\frac{(2n)!}{n!}\nonumber\\
&\geq& e^{-1}~2^{-n}\frac{\sqrt{4\pi n}~(2n)^{2n}~e^{-2n}}{\sqrt{2\pi n}~n^{n}~e^{-n}}=\sqrt{2}e^{-1}~\left(\frac{2}{e}\right)^{n}~~n^n .\label{stirling}
\end{eqnarray}
By \cite[Proposition 11.6.6]{mf-dm-13}, the probability of the following
event
\begin{equation}\label{event-control}
(Z/z)^2\leq \frac{e^2}{\sqrt{2}}~\left[\frac{1}{2}+\left(\delta+\sqrt{\delta}\right)\right]
\end{equation}
is greater than $1-e^{-\delta}$, for any $\delta\geq 0$.

The above estimate also implies that
$$
\EE\left(\exp{\left(t Z^2\right)}\right)\leq \frac{e}{\sqrt{2}}~\EE\left(\exp{\left(t(z) V^2\right)}\right)=\frac{e}{\sqrt{2}}~\frac{1}{\sqrt{1-2t(z)}}\quad\mbox{\rm with}\quad
t(z)=\frac{e}{2}~z^2t<1/2
$$
from which we check that
$$
t(z)\leq (1-\epsilon)/2\Rightarrow\begin{array}[t]{rcl}
\EE\left(\exp{\left(t Z^2\right)}\right)&\leq& \frac{e}{\sqrt{2}}~\exp{\left(t(z)\left[1+\frac{t(z)}{1-2t(z)}\right]\right)}\\
&\leq&
\frac{e}{\sqrt{2}}~\exp{\left(t(z)\left[1+\frac{1}{2}\left(\frac{1}{\epsilon}-1\right)\right]\right)} =
\frac{e}{\sqrt{2}}~\exp{\left(\frac{t(z)}{2}\left[1+\frac{1}{\epsilon}\right]\right)}.
\end{array}
$$
In summary we have
$$
t\leq (1-\epsilon)/(z^2e)\Rightarrow \EE\left(\exp{\left(t Z^2\right)}\right)\leq \frac{e}{\sqrt{2}}~\exp{\left(\frac{e}{4}~z^2t\left[1+\frac{1}{\epsilon}\right]\right)}.
$$
Choosing $t=(1-\epsilon)/(z^2e)$, we conclude that
$$
\forall \epsilon\in ]0,1]\qquad
\EE\left(\exp{\left[\frac{1-\epsilon}{e} \left(\frac{Z}{z}\right)^2\right]}\right)\leq \frac{e}{\sqrt{2}}~\exp{\left(\frac{1}{4}~
\frac{1-\epsilon^2}{\epsilon}\right)}.
$$
\subsection{Proof of Theorem~\ref{gronwall-2dim}}\label{proof-gronwall}
When $\Ua_t=\Va_t=0$ we have
$$
d\Vert\Xa_t\Vert^2\leq \rho(\Aa_t)~\Vert\Xa_t\Vert^2~dt+\sqrt{\Wa_t}~\Vert\Xa_t\Vert^2~d\Na_t\quad
$$
with the martingale
$$
d\Na_t:=\frac{1}{\sqrt{\Wa_t}~\Vert\Xa_t\Vert^{2}}~ d\Ma_t
~\Longrightarrow~
\partial_t\langle  \Na\rangle_t\leq 1 .
$$
This implies that
\begin{eqnarray*}
d\log{\Vert \Xa_t\Vert^2}&\leq & \Vert \Xa_t\Vert^{-2}~\left( \rho(\Aa_t)~\Vert\Xa_t\Vert^2~dt+\sqrt{\Wa_t}~\Vert\Xa_t\Vert^2~d\Na_t\right)-\frac{1}{2}~\Wa_t~dt\\
&=&\left(\rho(\Aa_t)-\frac{1}{2}~\Wa_t\right)~dt+\sqrt{\Wa_t}~d\Na_t
\end{eqnarray*}
from which we prove that
$$
\exp{\left(-\int_0^t\rho(\Aa_s)~ds\right)}~\Vert \Xa_t\Vert^2\leq \Vert \Xa_0\Vert^2~\Ea_t
$$
with the exponential martingale
$$
\Ea_t:=\exp{\left(\int_0^t\sqrt{\Wa_s}~d\Na_s-\frac{1}{2}~\int_0^t\Wa_s~ds\right)}.
$$

Next we provide a proof of the second assertion based on the above formula. 
For any $n\geq 0$, we observe that
$$
\exp{\left(-n\int_0^t\rho(\Aa_s)~ds\right)}~\Vert \Xa_t\Vert^{2n}\leq \exp{\left(\frac{n(n-1)}{2}~\int_0^t\Wa_s~ds\right)}~\Vert \Xa_0\Vert^{2n}~
\Ea_t(n)
$$
with the collection of exponential martingales
$$
\Ea_t(n):=\exp{\left(n\int_0^t\sqrt{\Wa_s}~d\Na_s-\frac{n^2}{2}~\int_0^t\Wa_s~ds\right)}.
$$
This implies that
$$
\EE\left(\exp{\left(-n~\int_0^t\left(\rho(\Aa_s)+\frac{(n-1)}{2}~\Wa_s\right)~ds\right)}~\Vert \Xa_t\Vert^{2n}~|~\Fa_0\right)\leq \Vert \Xa_0\Vert^{2n} .
$$
Arguing as above we use the decomposition
\begin{eqnarray*}
\EE\left(\Vert \Xa_t\Vert^n~|~\Fa_0\right)&=&
\EE\left(
\exp{\left(\frac{n}{2}~\int_0^t\left(\rho(\Aa_s)+\frac{(n-1)}{2}~\Wa_s\right)~ds\right)}~\right.\\
&&\left.\hskip3cm
\exp{\left(-\frac{n}{2}~\int_0^t\left(\rho(\Aa_s)+\frac{(n-1)}{2}~\Wa_s\right)~ds\right)}~\Vert \Xa_t\Vert^{n}~|~\Fa_0\right)
\end{eqnarray*}
to check that
$$
\EE\left(\Vert \Xa_t\Vert^n~|~\Fa_0\right)\leq 
\EE\left(
\exp{\left(n~\int_0^t\left(\rho(\Aa_s)+\frac{(n-1)}{2}~\Wa_s\right)~ds\right)}~
~|~\Fa_0\right)^{1/2}~\Vert  \Xa_0\Vert^{n} .
$$
This ends the proof of the first assertion.

More generally, we have
$$
d\Vert\Xa_t\Vert^2\leq \left[\rho(\Aa_t)~\Vert\Xa_t\Vert^2+\Ua_t\right]~dt+d\Ma_t .
$$
This yields
$$
\begin{array}{l}
d\Vert\Xa_t\Vert^{2n}\\
\\
\leq n~\Vert\Xa_t\Vert^{2(n-1)}~d\Vert\Xa_t\Vert^2+
\frac{n(n-1)}{2}~\Vert\Xa_t\Vert^{2(n-2)}~\left[ \Wa_t~\Vert\Xa_t\Vert^4+\Va_t~\Vert\Xa_t\Vert^2\right]~dt\\
\\
\leq -\Lambda_{n}(\Aa_t,\Wa_t)~\Vert\Xa_t\Vert^{2n}dt+n~\left[\frac{(n-1)}{2}~\Va_t~+\Ua_t~\right]~\Vert\Xa_t\Vert^{2(n-1)}~dt+n~\Vert\Xa_t\Vert^{2(n-1)}~d\Ma_t
\end{array}
$$
with
$$
-\Lambda_{n}(\Aa_t,\Wa_t):=~n~\rho(\Aa_t)~+\frac{n(n-1)}{2}~ \Wa_t .
$$
Observe that
\begin{eqnarray*}
\Lambda_{n}(\Aa_t,\Wa_t)-\Lambda_{n-1}(\Aa_t,\Wa_t)&=&-n~\rho(\Aa_t)~-n(n-1)~ \Wa_t/2
\\
&&\hskip2cm+(n-1)~\rho(\Aa_t)~+(n-1)(n-2)~ \Wa_t/2\\
&=&-\rho(\Aa_t)-(n-1)\Wa_t .
\end{eqnarray*}
This shows that
$$
\begin{array}{rcl}
\displaystyle\overline{\Ua}_t/\Ua_t&=&\displaystyle\exp{\left(\int_0^t\left[\Lambda_{n}(\Aa_s,\Wa_s)-\Lambda_{n-1}(\Aa_s,\Wa_s)\right]ds\right)}
=\overline{\Va}_t/\Va_t .
\end{array}$$
We set
$$
\begin{array}{rcl}
\displaystyle\Ya_t^n&:=&\displaystyle\exp{\left(\int_0^t\Lambda_{n}(\Aa_s,\Wa_s)ds\right)}~\Vert\Xa_t\Vert^{2n} .
\end{array}$$
Notice that
$$
\begin{array}{l}
\displaystyle\exp{\left(\int_0^t\Lambda_{n}(\Aa_s,\Wa_s)ds\right)}~\Ua_t~\Vert\Xa_t\Vert^{2(n-1)}\\
\\
=
\displaystyle\Vert\Xa_t\Vert^{2(n-1)}~\exp{\left(\int_0^t\Lambda_{n-1}(\Aa_s,\Wa_s)ds\right)}\\
\\
\hskip4cm
\displaystyle\times
\exp{\left(\int_0^t\left[\Lambda_{n}(\Aa_s,\Wa_s)-\Lambda_{n-1}(\Aa_s,\Wa_s)\right]ds\right)}~\Ua_t
=\Ya_t^{n-1}\overline{\Ua}_t .
\end{array}
$$
This shows that
$$
d\Ya_t^n\leq n~\Ya_t^{n-1}~\left[~\overline{\Ua}_t+\frac{(n-1)}{2}~\overline{\Va}_t~\right]~dt+d\overline{\Ma}_t
$$
with the martingale
$$
d\overline{\Ma}_t:=n~\exp{\left(\int_0^t\Lambda_{n}(\Aa_s,\Wa_s)ds\right)}~\Vert\Xa_t\Vert^{2(n-1)}~d\Ma_t .
$$
 This implies that
$$
\partial_t\EE\left(\Ya_t^n~|~\Fa_0\right)\leq n~\EE\left(\Ya_t^{n-1}~\left[~\overline{\Ua}_t+\frac{(n-1)}{2}~\overline{\Va}_t~\right]~|~\Fa_0\right) .
$$
Using H\"older inequality we have
$$
\EE\left(~\overline{\Ua}_t~\Ya_t^{n-1}~|~\Fa_0\right)\leq \EE\left(\overline{\Ua}_t^n~|~\Fa_0\right)^{1/n}
\EE\left(\Ya_t^{n}~|~\Fa_0\right)^{1-1/n} .
$$
This yields the estimate
$$
\partial_t\EE\left(\Ya_t^n~|~\Fa_0\right)\leq n~\EE\left(\Ya_t^{n}~|~\Fa_0\right)^{1-1/n}~\left[~\EE\left(\overline{\Ua}_t^n~|~\Fa_0\right)^{1/n}+\frac{(n-1)}{2}~\EE\left(\overline{\Va}_t^n~|~\Fa_0\right)^{1/n}~\right]
$$
and therefore
\begin{eqnarray*}
\partial_t\EE\left(\Ya_t^n~|~\Fa_0\right)^{1/n}&=&\frac{1}{n}~\EE\left(\Ya_t^n~|~\Fa_0\right)^{-(1-1/n)}~
\partial_t\EE\left(\Ya_t^n~|~\Fa_0\right)\\&\leq& \EE\left(\overline{\Ua}_t^n~|~\Fa_0\right)^{1/n}+\frac{(n-1)}{2}~\EE\left(\overline{\Va}_t^n~|~\Fa_0\right)^{1/n}~.
\end{eqnarray*}
We conclude that
$$
\EE\left(\Ya_t^n~|~\Fa_0\right)^{1/n}\leq \EE\left(\Ya_0^n~|~\Fa_0\right)^{1/n}+\int_0^t~\left[
 \EE\left(\overline{\Ua}_s^n~|~\Fa_0\right)^{1/n}+\frac{(n-1)}{2}~\EE\left(\overline{\Va}_s^n~|~\Fa_0\right)^{1/n}\right]~ds .
$$
Using the decomposition
$$
\EE\left(\Vert\Xa_t\Vert^{n}~|~\Fa_0\right)
=
\EE\left(\exp{\left(-\frac{1}{2}\int_0^t\Lambda_{n}(\Aa_s,\Wa_s)ds\right)}~~\exp{\left(\frac{1}{2}\int_0^t\Lambda_{n}(\Aa_s,\Wa_s)ds\right)}~\Vert\Xa_t\Vert^{n}~|~\Fa_0\right)
$$
and Cauchy-Schwartz inequality we check that
$$
\EE\left(\Vert\Xa_t\Vert^{n}~|~\Fa_0\right)^{2/n}
\leq 
\EE\left(\exp{\left(-\int_0^t\Lambda_{n}(\Aa_s,\Wa_s)ds\right)}~|~\Fa_0\right)^{1/n}
\EE\left(\Ya_t^n~|~\Fa_0\right)^{1/n} .
$$
This implies that
$$
\begin{array}{l}
\EE\left(\Vert\Xa_t\Vert^{n}~|~\Fa_0\right)^{2/n}
\displaystyle\leq \EE\left(\exp{\left(-\int_0^t\Lambda_{n}(\Aa_s,\Wa_s)ds\right)}\right)^{1/n}\\
\\
\hskip4cm\times\displaystyle\left[
\Vert\Xa_0\Vert^{2}+\int_0^t~\left[
 \EE\left(\overline{\Ua}_s^n~|~\Fa_0\right)^{1/n}+\frac{(n-1)}{2}~\EE\left(\overline{\Va}_s^n~|~\Fa_0\right)^{1/n}\right]~ds
\right]
\end{array}.$$

This ends the proof of the theorem.\cqfd

\subsection{Proof of Proposition~\ref{prop-Laplace-Chi2}}\label{proof-prop-appendix}

By (\ref{hilbert-G}), for any $m\geq 1$ we have the uniform estimate
$$
\EE\left(\Vert\Xa_t\Vert^{2m}\right)^{1/m}\leq (u_t(a)+mv_t(a))
\Rightarrow
2~\EE\left(\Vert\Xa_t\Vert^{2m}\right)\leq~(2u_t(a))^m+(2v_t(a))^m~m^m .
$$
Using Stirling approximation (\ref{stirling}) 
we have 
$$
(2\gamma~v_t(a))^m
m^m\leq  \frac{e}{\sqrt{2}}~\left(e\gamma v_t(a)\right)^{m}~2^{-m}\frac{(2m)!}{m!}\Rightarrow 
\sum_{m\geq 0}(2v_t(a))^m~m^m\leq \frac{e}{\sqrt{2}}~\frac{1}{\sqrt{1-2ev_t(a)\gamma}}
$$
for any $\gamma<1/(2ev_t(a))$.
This yields
$$
2~\EE\left(e^{\gamma \Vert\Xa_t\Vert^{2}}\right)\leq e^{2\gamma u_t(a)}+\frac{e}{\sqrt{2}}~\frac{1}{\sqrt{1-2e\gamma v_t(a)}} .
$$
Choosing $\gamma=(1-\epsilon)/(2ev_t(a))$, with $\epsilon \in ]0,1]$ we find that
$$
\EE\left(\exp{\left[\frac{(1-\epsilon)}{e}~\frac{1}{2v_t(a)}~\Vert\Xa_t\Vert^{2}\right]}\right)\leq \frac{1}{2}~
e^{\frac{1-\epsilon}{e} \frac{u_t(a)}{v_t(a)}}+\frac{e}{2\sqrt{2}}~\frac{1}{\sqrt{\epsilon}}.$$
This ends the proof of (\ref{expo-chi2-unif}). Now we come to the proof of (\ref{expo-chi2}).

We have
$$
d\Vert\Xa_t\Vert^2\leq \left[-a~\Vert\Xa_t\Vert^2+\Ua_t\right]~dt+d\Ma_t .
$$
This implies that
$$
\int_0^t~\Vert\Xa_s\Vert^2~ds\leq \int_0^t~e^{-as}~\left[\int_0^s~e^{au}~\Ua_u~du\right]~ds+
\int_0^t~e^{-as}~\left[\int_0^s~e^{au}~d\Ma_u\right]~ds .
$$
On the other hand, by an integration by part we have
$$
a\displaystyle\int_0^t~e^{-a s}\left(
\int_0^s~e^{a u}~\Ua_u~du\right)ds
\displaystyle =\int_0^t\left(1-e^{-a (t-s)}\right)~\Ua_s~ds
$$
and
$$
a\displaystyle\int_0^t~e^{-a s}\left(
\int_0^s~e^{a u}~d\Ma_{u}\right)ds
\displaystyle =\int_0^t\left(1-e^{-a (t-s)}\right)~~d\Ma_{s} .
$$
This implies that  
$$
a~\int_0^t~
\Vert\Xa_s\Vert^2~ds\leq \int_0^t\Ua_s~ds+\overline{\Ma}^{(t)}_t 
$$
with the terminal state $\overline{\Ma}^{(t)}_t$
of the collection of martingales $\overline{\Ma}^{(t)}_u$ on $[0,t]$ defined by 
$$
\begin{array}{l}
\displaystyle\forall 0\leq u\leq t\qquad
\overline{\Ma}^{(t)}_u:=\int_0^u\left(1-e^{-a (t-s)}\right)~~d\Ma_{s}\\
\\
\displaystyle\Longrightarrow\quad \partial_u\langle
\overline{\Ma}^{(t)}
\rangle_u\leq  v~\Vert\Xa_u\Vert^2 .
\end{array}
$$
Therefore for any $\gamma\geq 0$ we have
$$
\gamma~\left[\left(a-\frac{\gamma}{2}~v\right)~\int_0^t~\Vert\Xa_s\Vert^2~ds-
\int_0^t\Ua_s~ds\right]\leq 
\gamma\overline{\Ma}^{(t)}_t-\frac{\gamma^2}{2}\langle \overline{\Ma}^{(t)}\rangle_t .$$

This implies that
$$
\EE\left(\exp{\left[\gamma~\left[\left(a-\frac{\gamma}{2}~v\right)~\int_0^t~\Vert\Xa_s\Vert^2~ds-
\int_0^t\Ua_s~ds\right]\right]}~|~\Fa_0\right)\leq 1.
$$
Using the decomposition
$$
\begin{array}{l}
\displaystyle\exp{\left[\frac{\gamma}{2}~\left(a-\frac{\gamma}{2}~v\right)~\int_0^t~\Vert\Xa_s\Vert^2~ds\right]}
\\
\\
=
\displaystyle\exp{\left[\frac{\gamma}{2}~\left[\left(a-\frac{\gamma}{2}~v\right)~\int_0^t~\Vert\Xa_s\Vert^2~ds-
\int_0^t\Ua_s~ds\right]\right]}\times
\exp{\left[\frac{\gamma}{2}~
\int_0^t\Ua_s~ds\right]},\
\end{array}
$$
replacing $\gamma/2$ by $\gamma$, by Cauchy-Schwarz inequality we find that
$$
\EE\left(\displaystyle\exp{\left[v~\gamma~\left(\frac{a}{v}-\gamma\right)~\int_0^t~\Vert\Xa_s\Vert^2~ds\right]}\right)
\leq \EE\left(\exp{\left[2\gamma
\int_0^t\Ua_s~ds\right]}\right)^{1/2}
$$
for any $\gamma\leq {a}/{v}$.
Observe that
$$
\gamma~\left(\frac{a}{v}-\gamma\right):=\frac{\alpha}{v}\leq c^2
\quad\mbox{\rm with}\quad
c=\frac{a}{v} .$$
We also have
$$
\gamma~\left(\frac{a}{v}-\gamma\right)=\frac{\alpha}{v}\Longleftrightarrow\gamma\in\left\{c/2-\sqrt{(c/2)^2-{\alpha}/{v}},c/2+\sqrt{(c/2)^2-{\alpha}/{v}}
\right\} .
$$
Choosing the smallest value we prove that
\begin{eqnarray*}
\EE\left(\displaystyle\exp{\left[v~\alpha/v~\int_0^t~\Vert\Xa_s\Vert^2~ds\right]}\right)^2
&\leq& \EE\left(\exp{\left[\left((a/v)-\sqrt{(a/v)^2-{4\alpha}/{v}}\right)
\int_0^t\Ua_s~ds\right]}\right)\\
&=& \EE\left(\exp{\left[\frac{4\alpha/v}{(a/v)+\sqrt{(a/v)^2-{4\alpha}/{v}}}~
\int_0^t\Ua_s~ds\right]}\right)
\end{eqnarray*}
for any $\beta=\alpha/v\leq a^2/(2v)^2$, or equivalently
\begin{eqnarray*}
\EE\left(\displaystyle\exp{\left[v~\beta~\int_0^t~\Vert\Xa_s\Vert^2~ds\right]}\right)
&\leq&\EE\left(\exp{\left[\frac{4\beta}{(a/v)+\sqrt{(a/v)^2-4\alpha/v}}~
\int_0^t\Ua_s~ds\right]}\right)^{1/2} .
\end{eqnarray*}

This ends the proof of the proposition.
\cqfd

\section{Proof of formulae (\ref{example-non-linear}), (\ref{U1U2U}) and (\ref{2nd-formula-Lipschitz})}\label{proof-formula-Vnl}

We start with the proof of (\ref{example-non-linear}).
To clarify the presentation, we write $\Qa$ instead of $\Qa_2$.
Let $\langle x,y\rangle_{\Qa}=\langle x,\Qa y\rangle$ and $\langle x,x\rangle_{\Qa}^{1/2}=\Vert x\Vert_{\Qa}$ be the inner product and the norm associated with the symmetric definite positive matrix $Q$.
Let $z(t)=x+t(y-x)$ be an interpolating path from $x$ to $y$, indexed by $t\in [0,1]$. Also let
\begin{eqnarray*}
\varphi(t)&:=&\langle u,\partial^2 V(z(t))u\rangle=\Vert x+t(y-x)\Vert_{\Qa}~\Vert u\Vert^2_{\Qa}+\frac{1}{\Vert x+t(y-x)\Vert_{\Qa}}~ \langle u,x+t(y-x)\rangle_{\Qa}^2.
\end{eqnarray*}
Taking the derivative w.r.t. $t$ we find that
\begin{eqnarray*}
\partial\varphi(t)
&=&\frac{\langle x+t(y-x),y-x\rangle_{\Qa}}{\Vert x+t(y-x)\Vert_{\Qa}}~\Vert u\Vert_{\Qa}^2+\frac{2}{\Vert x+t(y-x)\Vert_{\Qa}}~ \langle u,x+t(y-x)\rangle_{\Qa}~ \langle u,y-x\rangle_{\Qa}\\
&&\hskip4cm-\frac{ \langle u,x+t(y-x)\rangle^2_{\Qa}}{\Vert x+t(y-x)\Vert_{\Qa}^3}~\langle x+t(y-x),y-x\rangle_{\Qa}\\
&=&\frac{\langle x+t(y-x),y-x\rangle_{\Qa}}{\Vert x+t(y-x)\Vert_{\Qa}}~\Vert u\Vert_{\Qa}^2~\left[1-\frac{ \langle u,x+t(y-x)\rangle_{\Qa}^2}{\Vert u\Vert_{\Qa}^2~\Vert x+t(y-x)\Vert_{\Qa}^2}~\right]\\
&&\hskip5cm+2~\frac{ \langle u,x+t(y-x)\rangle_{\Qa}~}{\Vert x+t(y-x)\Vert_{\Qa}\Vert u\Vert_{\Qa}}~ \langle u,y-x\rangle_{\Qa}~\Vert u\Vert_{\Qa}
\end{eqnarray*}
On the other hand we have
$$
-1\leq \frac{ \langle u,x+t(y-x)\rangle_{\Qa}~}{\Vert x+t(y-x)\Vert_{\Qa}~\Vert u\Vert_{\Qa}}\leq 1
\Longrightarrow
\vert\partial\varphi(t)\vert\leq 4~\Vert y-x\Vert_{\Qa}~\Vert u\Vert_{\Qa}^2.
$$
This yields
$$
\vert \langle u,\left[\partial^2 V(x)-\partial^2 V(y)\right]u\rangle\vert\leq 4~\Vert y-x\Vert_{\Qa}~\Vert u\Vert_{\Qa}^2\leq 4~\Vert \Qa\Vert^{3/2}~\Vert y-x\Vert~\Vert u\Vert^2
$$
This ends the proof of (\ref{example-non-linear}). 

Now we come to the proof of (\ref{U1U2U}).

\begin{eqnarray*}
\partial_{x_i}\Va(x)&=&\partial\Ua_1(x_i)+\sum_{1\leq k\leq r_1,~k\not=i}~\left(\partial_{1}\Ua_2(x_i,x_k)+\partial_{2}\Ua_2(x_k,x_i)\right)\\
\partial_{x_i,x_j}\Va(x)&=&\left\{
\begin{array}{lcl}
\partial_{2,1}\Ua_2(x_i,x_j)+\partial_{1,2}\Ua_2(x_j,x_i)
&\mbox{\rm for}& i\not=j\\
&&\\
\displaystyle\partial^2\Ua_1(x_i)+\sum_{1\leq j\leq r_1,~j\not=i}~\left(\partial_{2,2}+\partial_{1,1}\right)\Ua_2(x_i,x_j)&\mbox{\rm for}& i=j
\end{array}
\right.
\end{eqnarray*}
In the above display, $\partial\Ua_1$ and $\partial^2\Ua_1$ stands for the first and second derivative of $\Ua_1$ on $\RR$; and
$\partial_{i}\Ua_{2}$ stands for the partial derivatives  of $\Ua_2:(x_1,x_2)\in\RR^2\mapsto \Ua_2(x_1,x_2)\in [0,\infty[ $ with respect to the $i$-th coordinate $x_i$; and  $\partial_{i,j}:=\partial_{i}\partial_j=\partial_{j}\partial_i$, with $i,j\in\{1,2\}$.
\begin{eqnarray*}
\langle  \partial^2\Va(x)~s,s\rangle&=&\sum_{1\leq i\leq r_1}~\langle \partial^2\Ua_1(x_i)~s_i,s_i\rangle+
\sum_{1\leq i\not=j\leq r_1}\langle \partial^2\Ua_2(x_i,x_j)~(s_i,s_j),(s_i,s_j)\rangle\\
&\geq & u_1~\sum_{1\leq i\leq r_1}~s_i^2+u_2~
\sum_{1\leq i\not=j\leq r_1}~(s_i^2+s^2_j)=v~\langle  s,s\rangle .
\end{eqnarray*}
This ends the proof of (\ref{U1U2U}).

We end this section with the proof of (\ref{2nd-formula-Lipschitz}).
We have
$$
\vert \langle (s_1,s_2),\left[\partial^2\Ua_2(x_1,x_2)-\partial^2\Ua_2(y_1,y_2)\right](s_1,s_2)\rangle\vert\leq \kappa_{\partial^2\Ua_2}~\Vert (x_1,x_2)-(y_1,y_2)\Vert~\Vert (s_1,s_2)\Vert^2
$$
from which we find that
\begin{eqnarray*}
\vert \langle s, \left[ \partial^2\Va(x)- \partial^2\Va(y)\right]~s\rangle\vert&\leq&
\kappa_{\partial^2\Ua_1}~\sum_{1\leq i\leq r_1}~\vert x_i-y_i\vert~s_i^2\\
&&+
\kappa_{\partial^2\Ua_2}~\sum_{1\leq i\not=j\leq r_1}~\sqrt{(x_i-y_i)^2+(x_j-y_j)^2}~\left(s_i^2+s_j^2\right) .
\end{eqnarray*}
Using the fact that
$$
~\sum_{1\leq i\leq r_1}~\vert x_i-y_i\vert~\frac{s_i^2}{~\sum_{1\leq j\leq r_1}~s_j^2}\leq ~\left[\sum_{1\leq i\leq r_1}~\vert x_i-y_i\vert^2~\frac{s_i^2}{~\sum_{1\leq j\leq r_1}~s_j^2}\right]^{1/2}~\leq \Vert x-y\Vert
$$
and
$$
\sum_{1\leq i\not=j\leq r_1}~\sqrt{(x_i-y_i)^2+(x_j-y_j)^2}~\frac{\left(s_i^2+s_j^2\right)}{\sum_{1\leq k\not=l\leq r_1}(s_k^2+s_l^2)}
\leq \sqrt{2(r_1-1)}~\Vert x-y\Vert ,
$$
we prove that
\begin{eqnarray*}
\vert \langle s, \left[ \partial^2\Va(x)- \partial^2\Va(y)\right]~s\rangle\vert&\leq&
\kappa_{\partial^2\Ua_1}~\Vert x-y\Vert~\Vert s\Vert^2+
\kappa_{\partial^2\Ua_2}~(r_1-1)~ \sqrt{2(r_1-1)}~\Vert x-y\Vert~\Vert s\Vert^2.
\end{eqnarray*}
This ends the proof of  (\ref{2nd-formula-Lipschitz}).
\cqfd

Pierre DEL MORAL\\
INRIA Bordeaux Research Center (France) \& UNSW School of Mathematics and Statistics  (Australia)\\
p.del-moral@unsw.edu.au
\medskip

Aline KURTZMANN\\
Universit\'e de Lorraine, Institut Elie Cartan de Lorraine \\ CNRS, Institut Elie Cartan de Lorraine, UMR 7502, Vandoeuvre-l\`es-Nancy, F-54506, France.\\
aline.kurtzmann@univ-lorraine.fr
\medskip

Julian TUGAUT\\
Univ Lyon, Universit\'e Jean Monnet, CNRS UMR 5208, Institut Camille Jordan\\
 Maison de l'Universit\'e, 10 rue Tr\'efilerie, CS 82301, 42023 Saint-Etienne Cedex 2, France\\
 tugaut@math.cnrs.fr
 
\end{document}